\documentstyle{article}
\renewcommand{\a }{\alpha }
\renewcommand{\b }{\beta }
\renewcommand{\d }{\delta }

\newcommand{\e }{\varepsilon }
\newcommand{\g }{\gamma}
\newcommand{\G }{\Gamma }
\renewcommand{\l }{\lambda }

\newenvironment{pf}{\noindent{\sc
Proof}.\enspace}{\rule{2mm}{2mm}\medskip}

\newtheorem{Theorem}{Theorem}[section]

\newtheorem{Lemma}[Theorem]{Lemma}

\newtheorem{remark}[Theorem]{Remark}
\newtheorem{remarks}[Theorem]{Remarks}
\newenvironment{Remark}{\begin{remark}
\rm}{\rule{2mm}{2mm}\end{remark}}

\def\re{{\rm I\! R}}

\begin{document}

\title{ {\sc Bifurcation results for semilinear elliptic problems in
$\re^{N}$} \thanks{This research was supported by MURST ``Variational Methods
and Nonlinear
    Differential Equations''} }

   \author { Marino Badiale\thanks{Dipartimento di
Matematica,
   Universit\`a di Torino. Via Carlo Alberto 10, 10123 Torino, Italy.
   email: badiale@dm.unito.it} $\;$ and Alessio
Pomponio\thanks{SISSA,via Beirut 2-4, 34014 Trieste, Italy.
   email: pomponio@sissa.it} }

\date{}

\parindent=8mm

\maketitle

\noindent{\sc Abstract:}
In this paper we obtain, for a semilinear elliptic problem in 
$\re^{N}$, families of solutions bifurcating
from the bottom of the spectrum of $-\Delta$. The problem is 
variational in nature and we apply a nonlinear
reduction method which allows us to search for solutions as critical 
points of suitable functionals defined on
finite-dimensional manifolds.

\

\section{Introduction and Main Results}\label{sec:int}
An interesting problem in bifurcation phenomena is to look for 
solutions bifurcating not from an
eigenvalue but from a point of the
continuous spectrum of the the linearized operator of the involved equation.
Typical examples of differential operators
with continuous spectrum are the Laplace or the Schr\"odinger
operators in all $\re^{N}$, and there are now many results on
bifurcation of solutions for
semilinear elliptic equations in $\re^{N}$, for example see \cite{S}, 
\cite{S1}, \cite{S2}, \cite{S3}, \cite{M}.
See also \cite{S4}, and the references therein, for the study of 
bifurcation into spectral gaps.
A.Ambrosetti and the first author have studied such kind of problems 
in \cite{AB2} and \cite{AB3},
obtaining several results on
bifurcation of solutions for a one-dimensional differential equation. In this
paper we pursue such a study, generalizing some of the results of 
\cite{AB2} to higher dimensions, and considering
also the case of a critical nonlinearity. In section 5 of this paper 
we also fill a gap in the proof of theorem 3.2 in
\cite{AB2}. We thank S. Kr\"omer, who pointed out this gap, for his 
remarks and for several useful discussions.
\par \noindent We consider the equation

\begin{equation}\label{eq:gene}
\left\{
          \begin{array}{rcl}
-\Delta \psi - \l\psi &=&a(x) |\psi |^{p-1}\psi + b(x)|\psi
|^{q-1}\psi ,\quad x\in \re^{N},\\
\lim_{|x|\to \infty}\psi(x)& =&0,
            \end{array}
\right.
\end{equation}

\noindent where $N\geq 1$, $\l$ is a negative parameter, $1<p<q\leq
\frac{N+2}{N-2}$ if $N\geq 3$
(and $q<+ \infty$ if $N=1, 2$), $p< 1+ 4/N$ and $a,b : \re^{N}
\rightarrow \re$ satisfy suitable hypotheses (see below). Equation 
(\ref{eq:gene}) is an homogeneous equation,
so $\psi
=0$ is a solution for all $\l$,
the line $\{ ( \l , \psi =0) \, \vert \, \l \in \re \}$ is a line of 
trivial solutions and, as $q>
p>1$, the linearized operator at
$\psi =0$ is given by $\psi \rightarrow- \Delta \psi - \l \psi $. It 
is well known that $[0,+\infty )$
is the spectrum of $-\Delta $ on $\re^{N}$, and that it contains no 
eigenvalue. We
will find solutions bifurcating
from the bottom of the essential spectrum of $-\Delta$. To be
precise, by ``solution'' we mean a couple
$(\l , \psi_{\l})$ such that $\psi_{\l} \in H^{1}( \re^{N} )$ and
$\psi_{\l}$ is a solution of (\ref{eq:gene})
in the weak sense of $H^{1}( \re^{N} )$. We look for solutions
bifurcating from the origin in $H^{1}( \re^{N} )$, that is families 
$(\l , \psi_{\l})$ of
solutions of (\ref{eq:gene}) such that $\l \in (\l_{0}, 0 )$ for some 
$\l_{0}<0$ and $ \psi_{\l}
\rightarrow 0$ in $
H^{1}( \re^{N} )$ as $\l \rightarrow 0$.
\par \noindent Now let us state
the hypotheses on the functions $a,b$.
On $a$ we assume that there is $A >0$ such that either $a-A \in
L^{1}(\re^{N})$ or $a-A$
is asymptotic, at infinity, to $1/ \vert x \vert^{\gamma}$, for
suitable $\g$. To be precise, in the first case
we assume the following set of hypotheses.

  \begin{description}
        \item[${\bf (a_1 )}$]  $a-A$ is continuous, bounded and $a(x)-A\in
L^1(\re^{N})$.

      \item[${\bf (a_2 )}$] $\int_{\re^{N}} (a(x)-A) dx \not= 0$.

\end{description}

\noindent In the second case we assume the following hypothesis:

  \begin{description}
        \item[${\bf (a_3 )}$]  $a-{A}$ is continuous and there exist $L\not=
0 $ and $\g \in \, ]0, N [$ such that
        $\vert x \vert^{\g} (a(x)-A) \rightarrow L $ as $\vert x \vert
\rightarrow +\infty$.
        \end{description}

\noindent Notice that ${\bf (a_1 )}$ of course implies that $a-A\in
L^{p}(\re^{N})$
for all $p \in \, [1, +\infty ]$, while ${\bf (a_3 )}$ implies that
$a-A$ is bounded. For
$b$ we use some of the following assumptions.

  \begin{description}
        \item[${\bf (b_{1} )}$]  $b$ is continuous and bounded.

       \item[${\bf (b_{2} )}$] $b\in L^{\frac{2N}{N+2}}(\re^{N})$. If
$N \geq 2\frac{q-p}{p-1}$ we also assume that there
     exists $\b \in [1, \b^{*} [ $ such that
     $b\in L^{\b}(\re^{N})$, where
     $$\b^{*}= \frac{N(p-1)}{N(p-1) -2(q-p)}  \quad \hbox{if}
\quad N > 2\frac{q-p}{p-1},
     \quad \quad \b^{*}= +\infty \quad \hbox{if} \quad N =
2\frac{q-p}{p-1}.$$

       \item[${\bf (b_{3} )}$] $b\in L^{\frac{2N}{N+2}}(\re^{N})$. If
$\gamma \geq 2\frac{q-p}{p-1}$ we also assume that there
     exists $\b \in [1, \b^{*} [ $ such that
     $b\in L^{\b}(\re^{N})$, where
     $$\b^{*}= \frac{N (p-1)}{ \g (p-1) -2(q-p)} \quad \quad
\hbox{if} \quad \g > 2\frac{q-p}{p-1},
     \quad \quad \b^{*}= +\infty \quad \quad \hbox{if} \quad \g =
2\frac{q-p}{p-1}.$$

\end{description}

\bigskip

\noindent The value $\g$ in ${\bf (b_{3} )}$ is that given in ${\bf
(a_{3} )}$. We will assume either
${\bf (b_{1} )}$ and ${\bf (b_{2} )}$, or ${\bf (b_{1} )}$ and ${\bf
(b_{3} )}$.
Notice that, assuming ${\bf (b_{1} )}$, hypotheses ${\bf (b_{2} )}$
and ${\bf (b_{3} )}$ are obviously satisfied when
$b\in L^{1}(\re^{N})$.

\bigskip

\noindent We can now state our main results.

\begin{Theorem}\label{th:1}
Assume $1<p<q\leq \frac{N+2}{N-2}$ if $N\geq 3$, and $q<+\infty$ if
$N=1,2$. Suppose that ${\bf (a_1 )}$, ${\bf (a_2 )}$, ${\bf (b_1 )}$, ${\bf
(b_2 )}$ hold. Then (\ref{eq:gene})
has a family of solutions
bifurcating  from the origin in $ L^{\infty}(\re^{N})$. If, besides, 
$p<1+ \frac{4}{N}$, this family of
solutions bifurcates from the origin also in $ H^{1}(\re^{N})$.
\end{Theorem}
\bigskip

\begin{Theorem}\label{th:2}
Assume $1<p<q\leq \frac{N+2}{N-2}$ if $N\geq 3$, and $q<+\infty$ if
$N=1,2$. Suppose that ${\bf (a_3 )}$, ${\bf (b_1 )}$ and ${\bf (b_3 )}$ hold.
Then (\ref{eq:gene})
has a family of solutions
bifurcating  from the origin in $ L^{\infty}(\re^{N})$. If, besides, 
$p<1+ \frac{4}{N}$, this family of
solutions bifurcates from the origin also in $ H^{1}(\re^{N})$.
\end{Theorem}

\bigskip

\begin{Remark}
When $p\geq 1+ 4/N$, in $H^{1}(\re^{N})$ the solutions can bifurcate 
from infinity or can be bounded away both
from zero and infinity.
\end{Remark}

\bigskip

\begin{Remark}
An interesting question is to know if the solutions that we find form
a curve. We give some results in this
direction in section 5.
\end{Remark}

\bigskip

\noindent In the proof of theorems \ref{th:1} and \ref{th:2} we
follow the framework of \cite{AB2},
concerning the existence of critical points of perturbed functionals.
We start by a change of variables.
Let us set $u(x)=\e^{2/(1-p)}\psi(x/\e)$, $\l=-\e^{2}$, so that
equation (\ref{eq:gene})
becomes
\begin{equation}\label{eq:gene1}
-\Delta u + u =A |u|^{p-1}u +( a(x/\e)-A)|u|^{p-1}u  +
\e^{2\frac{q-p}{p-1}}b(x) \vert u \vert^{q-1}u.
\end{equation}

\noindent It is obvious that to any family $u_\e \in H^{1}(\re^{N})$
of solutions of
(\ref{eq:gene1}), bounded as $\e \to 0$, there corresponds
a family
$\psi_\e (x)=\e^{2/(p-1)}u_\e(\e x)$ of solutions of (\ref{eq:gene}). 
When $p<1+ 4/N$ it is easy to check that
$\psi_\e (x) \rightarrow 0$ in
$H^{1}(\re^{N})$, as $\e \to 0$. When
$p\geq 1+ 4/N$ we still get solutions, and it is easy to see that 
they vanish, as $\e
\rightarrow 0$, in $L^{\infty}(\re^{N})$, but they do not vanish in 
$L^{2}(\re^{N})$.
  Throughout this paper we will look for bounded families of
$H^{1}$-solutions of
  (\ref{eq:gene1}).
\medskip

\noindent The paper is organized as follows: after the introduction
(section 1) we give in section 2
a brief sketch of the abstract critical point theory for perturbed
functionals that we use to prove
theorems \ref{th:1} and \ref{th:2}. In section 3 we prove theorem
\ref{th:1} and in section 4 we prove
theorem \ref{th:2}. In section 5 we give some results on the
existence
of curves of solutions bifurcating from $(0,0)$, and we fill a gap in 
the proof of theorem 3.2 in \cite{AB2}.

\bigskip

{\bf Notation}
\smallskip
\par We collect below a list of the main notation used
throughout the paper.
\smallskip

\begin{description}
\item{$\bullet$}
  If $E$ is a Banach space, $F:E \rightarrow E$, and $u\in E$, then 
$DF(u) : E \rightarrow E$,
$D^{2}F(u):  E\times E \rightarrow E$ and
$D^{3}F(u):  E\times E \times E \rightarrow E$ are the first, second 
and third differential of $F$ at $u$,
which are
respectively linear, bilinear and three-times linear.
\medskip

\item{$\bullet$}
$L(E,E)$ is the space of linear continuous operators
from $E$ to $E$.

\medskip

\item{$\bullet$}
   $2^{*} = \frac{2N}{N-2}$ is the critical exponent for the
Sobolev embedding, when $N\geq 3$.
\medskip

\item{$\bullet$}
  We will use $C$ to denote any positive constant, that can
change from line to line.
\end{description}

\bigskip

\section{Abstract theory for perturbed functionals}

In this section we give the main ideas and results of a variational
method to study
critical points of perturbed functionals. The method has been
developed in
\cite{ABC}, \cite{AB1} , \cite{AB2} and then has been applied
to many different problems, see \cite{AGP1}, \cite{AGP2}, 
\cite{AMN1}, \cite{AMS}, \cite{BB}, \cite{BD},
\cite{BPV}.
We deal with a family of functionals $f_{\e}$, defined on a
  Hilbert space $E$, of the form
$$
f_{\epsilon}(u)=\frac{1}{2}\|u\|^2- F(u) + G(\e,u),
$$
where $\|\cdot \|$ is the norm in $E$, $F:E\to \re$ and $G:\re\times
E\to
\re$. We need the following hypotheses
\begin{description}
        \item[${\bf (F_0)}$] $F\in C^2$;

        \item[${\bf (G_0)}$]  $G$ is continuous in $(\e,u)\in \re\times E$
and $G(0,u)=0$ for all $u\in E$;

        \item[${\bf (G_1)}$]
$G$ is of class $C^2$ with respect to $u\in E$.

\end{description}

\

\noindent We will use the notation $F'(u)$,
respectively $G'(\e,u)$, to denote the functions defined by setting
$$
(F'(u)|v)=DF(u)[v],\quad \forall \;v\in E,
$$
and, respectively,
$$
(G'(\e,u)|v)=D_u G(\e,u)[v],\quad \forall \;v\in E,
$$
\noindent  where $(\cdot \, \vert \, \cdot )$ is the scalar product 
in $E$. Similarly, $F''(u)$, resp.
$G''(\e,u)$, denote
the maps in $L(E,E)$ defined by
$$
(F''(u)v | w) = D^2 F(u)[v,w] \quad \quad
(G''(\e,u)v|w)  =  D^2_{u u} G(\e,u)[v,w].
$$

\noindent In section 5 we will assume that $F,G$ are $C^{3}$. In this 
case we will denote $F'''(u)$,
$G'''(\e,u)$ the
bilinear maps defined by

$$ (F'''(u)[v_{1}, v_{2}] \, \vert \, v_{3} ) = D^{3}F(u) [v_{1}, 
v_{2}, v_{3}],
\quad \quad
(G'''(\e ,u)[v_{1}, v_{2}] \, \vert \, v_{3} ) = D^{3}_{uuu}G(\e , u) 
[v_{1}, v_{2}, v_{3}].$$

\noindent We also assume that $F$ satisfies

\

\begin{description}
        \item [${\bf (F_1 )}$]  there exists a $d$-dimensional $C^2$
manifold $Z$,
$d\geq 1$, consisting of critical       points of $f_0$, namely such that

$$\label{F1}
z-F'(z)=0,\quad \forall\; z\in Z.$$

\end{description}

\

\noindent Such a $Z$ will be called a {\sl critical manifold} of
$f_0$.

Let $T_z Z$ denote the tangent space to $Z$ at $z$ and $I_E$ denote
the Identity
map in $E$. We further suppose

\

\begin{description}
        \item[${\bf (F_2)}$]   $F'' (z)$ is compact $\forall \;z\in Z$;

        \item[${\bf (F_3)}$]  $T_z Z = Ker [I_E-F '' (z)]$, $\forall \;z\in
Z$.
\end{description}

\noindent We make the following further assumptions on $G$.

\

\begin{description}

        \item[${\bf (G_2)}$]   The maps $(\e,u) \mapsto G'(\e,u)$,
$(\e,u) \mapsto G''(\e,u)$
are continuous  (as maps from $\re\times E$ to $E$, resp. to
$L(E,E)$).
\end{description}

\

\begin{description}
\item[${\bf (G_3)}$]
there exist $\a>0$ and a continuous function $\G:Z\to \re$ such
that, for all $z \in Z$,
$$
\G(z)=\lim_{\e\to 0} \frac{G(\e,z)}{\e^\a}.
$$
and
$$
G'(\e,z)=o(\e^{\a/2}).
$$
\end{description}

\

\noindent In \cite{AB2} (see also \cite{AB1}, \cite{ABC}) the
following theorem is proved.

  \begin{Theorem}\label{th:abstr}
Suppose ${\bf (F_0 - F_3)}$ and ${\bf (G_0 - G_3)}$ hold and
assume there exist $\d>0$ and $z^*\in Z$ such that

\begin{equation}\label{eq:nondeg}
\quad\mbox{either}\; \; \min_{\|z-z^*\|=\d} \G(z) >
\G(z^*),\quad\mbox{or}\;\; \max_{\|z-z^*\|=\d} \G(z)\;<\G(z^*).
\end{equation}

\noindent Then, for $\e$ small, $f_\e$ has a critical point $u_\e $.
\end{Theorem}
\begin{pf}
We give only a sketch of the proof, divided in three steps.

\medskip
{\bf Step 1. } Using the Implicit Function Theorem one can find
$w=w(\e,z) \perp T_z Z $ such that
\begin{equation}\label{eq:w}
f^{'}_\e(z+w)\in T_z Z    , \quad  \vert
\vert w \vert \vert
=o (\e^{\a /2} )\quad \mbox{and} \quad \vert \vert D_{z}w( \e , z ) 
\vert \vert \rightarrow 0
\quad \mbox{as} \quad \e \rightarrow 0.
\end{equation}
\noindent Letting $Z_\e=\{z+w(\e,z)\}$, it turns out that
$Z_\e$ is locally diffeomorphic to $Z$ and any critical point of
$f_\e$ restricted to
$Z_\e$ is a stationary point of $f_\e$.

\medskip
{\bf Step 2. } Using the Taylor expansion we obtain, for $u= z+
w(\e , z) \in Z_{\e}$,

$$ f_{\e}(u) = c+ \e^{\a} \Gamma (z) + o(\e^{\a}),$$

\noindent where $c$ is a constant.

\medskip
{\bf Step 3. } It readily follows that, for small $\e$'s, $f_{\e}$
has a local
constrained minimum (or maximum) on $Z_{\e}$
at some $u_{\e}= z_{\e}+ w(\e, z_{\e}) \in Z_{\e}$, with $\vert \vert 
z_{\e}- z^{*} \vert \vert < \delta$.
According to step 1, such
$u_{\e}$ is a critical
point of $f_{\e}$ .

\end{pf}

\section{First bifurcation result.}

In this section we prove theorem \ref{th:1}. We want to apply the
abstract tools of the previous section,
and we start to set

$$E= H^{1}(\re^{N}) , \quad \quad \vert \vert u \vert \vert^{2} =
\int_{\re^{N}} \left( \vert \nabla u \vert^{2} +u^{2} \right) dx ,
\quad \quad
F(u) = \frac{A}{p+1} \int_{\re^{N}} \vert u \vert^{p+1} dx ,$$

\noindent and $G= G_{1 } +G_{2}$ where

$$
G_{1}(\e,u)=\left\{
          \begin{array}{lcl}
\frac{-1}{p+1}\int_{\re^{N}}( a(x/\e)-A)|u|^{p+1}dx &\mbox{if}& \e
\not= 0,\\
0&\mbox{if}& \e = 0
            \end{array}
\right.
$$

\noindent and
$$
G_{2}(\e,u)=\left\{
          \begin{array}{lcl}
\frac{-1}{p+1}\,
\e^{2\frac{q-p}{p-1}}\int_{\re^{N}}b(x/\e)|u|^{q+1}dx
&\mbox{if}& \e \not= 0,\\
0&\mbox{if}& \e = 0
            \end{array}
\right.
$$

\bigskip

\noindent Throughout this section we assume $N\geq 3$ and, of course,
$1<p<q\leq \frac{N+2}{N-2}$. The cases $N=1,2$ can
be handled in the same way, and in fact are easier.
We have now to verify that the hypotheses ${\bf (F_{0}-F_{3})}$ and
${\bf (G_{0}-G_{3})}$ are satisfied. The fact
that $q>p>1$ gives of course ${\bf (F_{0})}$ and ${\bf (G_{1})}$. It
is also well known (see \cite{BL1}, \cite{BL2}, \cite{K}) that
there exists a unique positive radial solution $z_{0}$ of
$$ - \Delta u + u = A \vert u \vert^{p-1}u , \quad \quad x \in
\re^{N} ,$$

\noindent  that $z_{0}$ is strictly radial decreasing, has an
exponential decay at infinity together with its derivatives, and that $f_{0} $
possesses a $N-$dimensional manifold of  critical points
$$Z= \{ z_{\theta}(x) = z_{0 }(x+ \theta ) \, \vert \,  \theta \in
\re^{N}  \, \} .$$

\noindent Furthermore, we know (see \cite{Oh}, \cite{ABC} and the references
therein) that $T_{z_{\theta}}Z = $ker$(I_{E}-
F'' (z_{\theta}))$ for all $z_{\theta} \in Z$. It is also easy
to check that $F'' (z_{\theta})$ is compact, for all $z_{\theta} \in 
Z$. In this way
all the hypotheses on $F$ are satisfied, and
the rest of this section is devoted to prove those on $G$. We will
get this by several lemmas. Let us prove
as first thing that the hypothesis ${\bf (G_{0})}$ is satisfied.
\bigskip

\begin{Lemma}\label{le:1}
Assume ${\bf (a_{1} ) }$ and ${\bf (b_{1} ) }$. Then $G$ is
continuous.
\end{Lemma}
\begin{pf}
We prove first that $G_{1}$ is continuous.
Assume that $(\e , u ) \rightarrow (\e_{0} , u_{0}  )$ in $\re \times
H^{1} (\re^{N} )$, with
$\e_{0}\not= 0$. Then we can write

$$ (p+1) \vert  G_{1} (\e , u ) - G_{1}(\e_{0} , u_{0}  ) \vert  =$$

$$ \Big\vert  \int_{\re^{N}} \left(   a\left(   \frac{x}{\e} \right)
-A  \right)  \, \vert u \vert^{p+1} \,
dx -  \int_{\re^{N}}  \left(   a\left(   \frac{x}{\e_{0} } \right)
-A  \right)  \, \vert u_{0}
\vert^{p+1} \, dx
  \Big\vert  \leq $$

$$ \Big\vert  \int_{\re^{N}}  \left(   a\left(   \frac{x}{\e} \right)
-A  \right)  \, \vert u \vert^{p+1} \,
dx - \int_{\re^{N}} \left(   a\left(   \frac{x}{\e_{0} } \right) -A
\right)  \, \vert u \vert^{p+1} \,
  dx \Big\vert  + $$

$$ \Big\vert  \int_{\re^{N}} \left(   a\left(   \frac{x}{\e_{0} }
\right) -A  \right)  \, \vert u \vert^{p+1} \,
dx -  \int_{\re^{N}}  \left(   a\left(   \frac{x}{\e_{0} } \right)
-A  \right)  \, \vert u_{0}
\vert^{p+1} \, dx \Big\vert  \leq $$

$$ \int_{\re^{N}} \Big\vert a\left(   \frac{x}{\e} \right) -
a\left(   \frac{x}{\e_{0} } \right)
\Big\vert  \, \vert u \vert^{p+1} \, dx +$$

$$C \int_{\re^{N}} \Big\vert  \, \vert u \vert^{p+1} - \vert u_{0}
\vert^{p+1} \Big\vert \, dx ,$$

\noindent By hypothesis $a$ is continuous and bounded, so it is easy
to deduce, by dominated convergence, that
the first term goes to zero, while the second one goes to zero by
hypothesis. Hence we deduce $ \vert
G_{1} (\e , u ) - G_{1}(\e_{0} , u_{0}  ) \vert  \rightarrow 0 .$
\par \noindent Now assume that $(\e , u ) \rightarrow (0, u_{0}  )$.
By definition $G_{1}(0, u_{0} )= 0$
and we have, applying H\"older inequality,

$$ (p+1) \vert G_{1} (\e , u ) \vert  \leq \int_{\re^{N}} \Big\vert
a\left(   \frac{x}{\e} \right) - A
\Big\vert  \, \vert u \vert^{p+1} \, dx \leq $$

$$ \left(\int_{\re^{N}} \Big\vert a\left(   \frac{x}{\e} \right) - A
\Big\vert^{\frac{2^*}{2^* - p-1} }dx \right)^{\frac{2^* -p-1}{2^* }
}\,
\left(\int_{\re^{N}} \vert u \vert^{2^*}dx \right)^{\frac{p+1}{2^*}}
.$$

\noindent By the change of variables $y= x/ \e $ we get

$$  \left(\int_{\re^{N}} \Big\vert a\left(   \frac{x}{\e} \right) - A
\Big\vert^{\frac{2^*}{2^* - p-1} }dx \right)^{\frac{2^* -p-1}{2^* } }=
\e^{N \frac{2^* -p-1}{2^*}} \,  \left(\int_{\re^{N}} \vert a (y) - A
\vert^{\frac{2^*}{2^* - p-1} }dy \right)^{\frac{2^* -p-1}{2^* } }.$$

\noindent As $a-A \in L^{\frac{2^* }{2^* -p-1} } (\re^{N})$ and $u
\in
H^{1}(\re^{N}) \subset
L^{2^* } (\re^{N}),$ we get $G_{1}(\e , u)  \rightarrow  0$ as $(\e ,
u)  \rightarrow ( 0, u_{0})$.

\par \noindent As to $G_{2}$, we argue in the same way. If $(\e , u )
\rightarrow (\e_{0} , u_{0}  )$ with
$\e_{0}\not= 0$, then

$$ (q+1) \vert  G_{2} (\e , u ) - G_{2}(\e_{0} , u_{0}  ) \vert  =$$

$$ \Big\vert  \e^{2\frac{q-p}{p-1}}   \int_{\re^{N}}   b\left(
\frac{x}{\e} \right)   \, \vert u \vert^{q+1} \,
dx -  \e_{0}^{2\frac{q-p}{p-1}}   \int_{\re^{N}}   b\left(
\frac{x}{\e_{0} } \right)  \, \vert u_{0}
\vert^{q+1} \, dx
  \Big\vert  \leq $$

$$\e^{2\frac{q-p}{p-1}}   \Big\vert  \int_{\re^{N}}    b\left(
\frac{x}{\e} \right)  \, \vert u \vert^{q+1} \,
dx - \int_{\re^{N}}   b\left(   \frac{x}{\e_{0} } \right)  \, \vert
u_{0}  \vert^{q+1} \, dx \Big\vert  + $$

$$ \vert \e^{2\frac{q-p}{p-1}} -\e_{0}^{2\frac{q-p}{p-1}} \vert \, \,
\Big\vert  \int_{\re^{N}}  b\left(   \frac{x}{\e_{0} } \right)  \,
\vert u_{0} \vert^{q+1} \, dx \Big\vert . $$

\noindent The first term can be treated as above, while the second
one obviously goes to zero as $\e \rightarrow
  \e_{0}$.
\par \noindent If $(\e , u ) \rightarrow (0, u_{0}  )$, we have

$$ (q+1) \vert  G_{2} (\e , u ) \vert \leq \e^{2\frac{q-p}{p-1}}
\int_{\re^{N}} \big\vert b\left(   \frac{x}{\e } \right) \big\vert \, \vert u
\vert^{q+1} \, dx \leq
C \e^{2\frac{q-p}{p-1}} .$$

\noindent So also $G_{2}$ is a continuous function, hence $G$ is
continuous and the lemma is proved.

\end{pf}

\bigskip

\noindent In the next lemma we prove that ${\bf (G_{2}) }$ is
satisfied.

\bigskip

\begin{Lemma} \label{le:2}
Assume ${\bf (a_{1} ) }$ and ${\bf (b_{1} ) }$.
Then $G'$ and $G''$ are continuous.
\end{Lemma}
\begin{pf}
Let us consider $G'_{1}$, and assume $(\e , u ) \rightarrow  (\e_{0}
, u_{0}  ) $ with $\e_{0} \not= 0$. We obtain
$$ \vert \vert  G'_{1}(\e , u )   - G'_{1}(\e_{0}  , u_{0}  )
\vert  \vert   =
\sup_{\vert \vert v \vert \vert \leq 1 }  \left\{ \big\vert  ((G'_{1}(\e , u
)   - G'_{1}(\e_{0}  , u_{0}  )  ) \vert
v )    \big\vert \right\} =$$

$$ \sup_{\vert \vert v \vert \vert \leq 1 }
\left\{  \Big\vert  \int_{\re^{N}} \left(   a\left(   \frac{x}{\e}
\right)
-A  \right)  \, \vert u \vert^{p-1}u \, v  \,
  dx - \int_{\re^{N}}  \left(   a\left(   \frac{x}{\e_{0} } \right)
-A  \right)  \, \vert u_{0}
  \vert^{p-1} u_{0} \, v  \, dx
  \Big\vert  \right\} \leq $$

$$ \sup_{\vert \vert v \vert \vert \leq 1 }
  \left\{ \Big\vert  \int_{\re^{N}} \left(   a\left(   \frac{x}{\e}
\right)
-A  \right)  \, \vert u \vert^{p-1}u \, v  \, dx
  - \int_{\re^{N}}  \left(   a\left(   \frac{x}{\e_{0} } \right) -A
\right)  \, \vert u
   \vert^{p-1} u\, v  \, dx
  \Big\vert  \right\} + $$

$$ \sup_{\vert \vert v \vert \vert \leq 1 }
\left\{ \Big\vert  \int_{\re^{N}} \left(   a\left(   \frac{x}{\e_{0}}
\right) -A  \right)  \, \vert u  \vert^{p-1}u
   \, v  \, dx - \int_{\re^{N}}  \left(   a\left(   \frac{x}{\e_{0} }
\right) -A  \right)  \, \vert u_{0}
    \vert^{p-1} u_{0} \, v  \, dx
  \Big\vert \right\} . $$

\noindent For the first term we can write

$$\Big\vert  \int_{\re^{N}} \left(   a\left(   \frac{x}{\e} \right)
-A  \right)  \, \vert u \vert^{p-1}u \, v  \,
  dx - \int_{\re^{N}}  \left(   a\left(   \frac{x}{\e_{0} } \right)
-A  \right)  \, \vert u
  \vert^{p-1} u\, v  \, dx
  \Big\vert  \leq$$

$$  \left(  \int_{\re^{N}} \Big\vert a\left(   \frac{x}{\e} \right) -
a\left(   \frac{x}{\e_{0} } \right)
\Big\vert^{\frac{p+1}{p}} \, \vert u \vert^{p+1}
\right)^{\frac{p}{p+1}}
\left(\int_{\re^{N}} \vert v \vert^{p+1} dx
\right)^{\frac{1}{p+1}} \leq$$

$$  C\, \left(  \int_{\re^{N}} \Big\vert a\left(   \frac{x}{\e}
\right) - a\left(   \frac{x}{\e_{0} } \right)
\Big\vert^{\frac{p+1}{p}} \, \vert u \vert^{p+1}
\right)^{\frac{p}{p+1}} ,$$

\noindent where $C$ is independent of $v$, ($\vert \vert v \vert
\vert \leq 1 $).
As above, this term tends to zero, by dominated
convergence. For the second term we have

$$  \Big\vert  \int_{\re^{N}} \left(   a\left(   \frac{x}{\e_{0}}
\right) -A  \right)  \, \vert u
\vert^{p-1}u  \, v  \, dx - \int_{\re^{N}}  \left(   a\left(
\frac{x}{\e_{0} } \right) -A
\right)  \, \vert u_{0}  \vert^{p-1} u_{0} \, v  \, dx
  \Big\vert  \leq$$

$$  C \, \int_{\re^{N}} \big\vert  \vert u \vert^{p-1}u - \vert
  u_{0}  \vert^{p-1}u_{0} \big\vert
\, \vert v \vert dx \leq
  C \, \left(  \int_{\re^{N}} \big\vert  \vert u \vert^{p-1}u -
  \vert u_{0}  \vert^{p-1}u_{0}
\big\vert^{\frac{p+1}{p}  } dx \right) ^{\frac{p}{p+1}  }
\left( \int_{\re^{N}}\vert v \vert^{p+1 } dx \right)^{\frac{1}{p+1}
} \leq$$

$$
C \, \left(  \int_{\re^{N}} \big\vert  \vert u \vert^{p-1}u - \vert
  u_{0}  \vert^{p-1}u_{0}
\big\vert^{\frac{p+1}{p}  } dx \right) ^{\frac{p}{p+1}  }, $$

\noindent and this term vanishes as $u \rightarrow u_{0} $. Hence we conclude
$\vert \vert  G_{1}'(\e , u) -   G_{1}'(\e_{0} , u_{0} )   \vert 
\vert \rightarrow 0$
as $(\e , u)\rightarrow (\e_{0} , u_{0} )$, $\e_{0} \not= 0$. Let us
now assume
$(\e , u ) \rightarrow  (0 , u_{0} ) $. By definition, $G_{1}' (0, u
) = 0$ and

$$\vert \vert G'_{1} (\e , u ) \vert \vert  =
\sup_{\vert \vert v \vert \vert \leq 1 }  \left\{ \Big\vert  \int_{\re^{N}}
\left(   a\left(   \frac{x}{\e} \right) -A
\right)  \, \vert u \vert^{p-1}u \, v  \, dx \Big\vert  \right\} \leq $$

$$  \sup_{\vert \vert v \vert \vert \leq 1 } \left\{ \int_{\re^{N}}
\left(     \Big\vert
a\left(   \frac{x}{\e} \right) -A  \Big\vert^{\frac{2^*}{2^* -p -1}}
dx \right)^{\frac{2^* -p-1}{2^*}}
\, \left(  \int_{\re^{N}} \vert u \vert^{2^*}dx \right)^{\frac{p}{2^*
}} \,
\left(  \int_{\re^{N}} \vert v \vert^{2^*}dx \right)^{\frac{1}{2^*
}}   \right\}  \leq $$

$$ C \, \e^{ N \frac{2^* -p-1}{2^*}} \left(  \int_{\re^{N}} \vert
a(y) - A)
\vert^{\frac{2^* }{2^* -p-1}} \, dy \right)^{\frac{2^* -p-1}{2^* }}\,
\left(  \int_{\re^{N}} \vert u \vert^{2^*}dx \right)^{\frac{p}{2^*
}}$$

\noindent and this term vanishes as $\e \rightarrow 0$. In this way
we have proved that $G_{1}'$ is continuous.
\par \noindent Similar arguments work for $G_{2}'$. Indeed, if $(\e ,
u )
\rightarrow  (\e_{0} , u_{0}  ) $ with $\e_{0} \not= 0$,
  we obtain

$$ \vert \vert  G'_{2}(\e , u )   - G'_{2}(\e_{0}  , u_{0}  )
\vert  \vert   =
\sup_{\vert \vert v \vert \vert \leq 1 } \big\vert ( (G'_{2}(\e , u
)   - G'_{2}(\e_{0}  , u_{0}  )  ) \vert
v ) \big\vert =$$

$$ \sup_{\vert \vert v \vert \vert \leq 1 } \left\{
  \Big\vert \e^{2   \frac{q-p}{p-1}} \int_{\re^{N}}   b\left(
\frac{x}{\e} \right)   \, \vert u \vert^{p-1}u \, v
   \, dx -  \e_{0}^{2   \frac{q-p}{p-1}}   \int_{\re^{N}}   b\left(
\frac{x}{\e_{0} }  \right)\, \vert u_{0}
    \vert^{p-1} u_{0} \, v  \, dx
  \Big\vert  \right\}  \leq $$

$$ \sup_{\vert \vert v \vert \vert \leq 1 }
  \left\{   \e^{2   \frac{q-p}{p-1}}  \Big\vert \int_{\re^{N}}
b\left(
\frac{x}{\e} \right)   \, \vert u \vert^{p-1}u \,
   v  \, dx -  \int_{\re^{N}}   b\left(   \frac{x}{\e_{0} } \right)
\, \vert u_{0}  \vert^{p-1} u_{0} \, v  \, dx
      \Big\vert \right\} + $$

$$ \sup_{\vert \vert v \vert \vert \leq 1 }
  \left\{ \Big\vert \e^{2   \frac{q-p}{p-1}} \int_{\re^{N}}
b\left(
\frac{x}{\e_{0}} \right)   \, \vert u_{0}
  \vert^{p-1}u_{0}  \, v  \, dx -  \e_{0}^{2   \frac{q-p}{p-1}}
\int_{\re^{N}}   b\left(   \frac{x}{\e_{0} }
   \right)    \, \vert u_{0}  \vert^{p-1} u_{0} \, v  \, dx
  \Big\vert \right\} \leq $$

$$  \e^{2   \frac{q-p}{p-1}}  \sup_{\vert \vert v \vert \vert \leq 1
}
  \left\{ \Big\vert \int_{\re^{N}}   b\left(   \frac{x}{\e} \right)
\, \vert
u \vert^{p-1}u \, v  \, dx -  \int_{\re^{N}}
   b\left(   \frac{x}{\e_{0} } \right)    \, \vert u_{0}  \vert^{p-1}
u_{0} \, v  \, dx   \Big\vert \right\} + $$

$$  \big\vert  \e^{2   \frac{q-p}{p-1}} -  \e_{0}^{2
\frac{q-p}{p-1}} \big\vert  \,
\sup_{\vert \vert v \vert \vert \leq 1 }  \int_{\re^{N}} \Big\vert
b\left(   \frac{x}{\e_{0} } \right) \,
  \vert u_{0}  \vert^{p-1} u_{0} \, v  \, \Big\vert dx .$$

\noindent The first term can be treated exactly as before, the second
term obviously vanishes as
$\e \rightarrow \e_{0} $. Let us now assume $(\e , u )\rightarrow (0,
u_{0})$. We obtain

$$  \vert \vert G_{2}' (\e ,u ) \vert \vert = \sup_{\vert \vert v \vert \vert
\leq 1} \big\vert \left( G'_{2} (\e, u)
\, \vert \, v \right) \big\vert \leq $$
$$  \sup_{\vert \vert v \vert \vert \leq 1} \e^{2 \frac{q-p}{p-1}}
\int_{\re^{N}}\big\vert
b\left( \frac{x}{\e} \right)\big\vert \,  \vert u \vert^{p} \vert v
\vert dx \leq C \, \e^{2 \frac{q-p}{p-1}} .$$

\par\noindent Now we have proved that $G'$ is continuous. The
argument to prove the continuity of $G''$
is almost the same and we leave it to the reader.
\end{pf}

\bigskip

\noindent Let us now verify that ${\bf (G_{3} )} $ is satisfied.

\bigskip

\begin{Lemma}\label{le:3}
Let us assume ${\bf  (a_{1}) }$, $({\bf  b_{1}
})$ and $({\bf  b_{2} })$. Let us define,
for $\theta \in \re^{N}$,

\begin{equation}\label{eq:gamdef}
  \Gamma (\theta ) = - \frac{1}{p+1} \, z_{0}^{p+1} (\theta )\,
\int_{\re^{N}} (a(y)-A )dy .
\end{equation}

\noindent Then

\begin{equation}\label{eq:lim1}
\lim_{\e \rightarrow 0} \frac{ G(\e , z_{\theta}) }{ \e^{N} } =
\Gamma (\theta )
\end{equation}

\noindent and

\begin{equation}\label{eq:lim2}
G'(\e , z_{\theta}) = O(\e^{\frac{N}{2}+1} ) .
\end{equation}

\end{Lemma}

\begin{pf}
As above we will study separately $G_{1}$ and $G_{2}$. By the change
of variables $y= \frac{x}{\e }$ we have

$$ G_{1}(\e , z_{\theta } )= - \frac{1}{p+1} \int_{\re^{N}}\left(
a\left( \frac{x}{\e }  \right)
-A \right) z_{0}^{p+1}(x+\theta ) dx= $$

$$- \frac{\e^{N}}{p+1}\int_{\re^{N}} (a(y)-A ) \, z_{0}^{p+1}(\e
y+\theta ) \, dy .$$

\noindent Since $a-A \in  L^{1} (\re^{N})$ and $z_{0}$ is bounded and
continuous, by dominated convergence we get
\begin{equation}\label{eq:lim3}
\lim_{\e \rightarrow 0} \frac{ G_{1}(\e , z_{\theta}) }{ \e^{N} } =
\Gamma (\theta )
\end{equation}

\noindent Hence, to prove (\ref{eq:lim1}) we have to show that

\begin{equation}\label{eq:lim4}
\lim_{\e \rightarrow 0} \frac{ G_{2}(\e , z_{\theta}) }{ \e^{N} } = 0.
\end{equation}

\noindent We distinguish two cases. Assume first $N< 2
\frac{q-p}{p-1} $. In this case

$$ \e^{-N} G_{2}(\e , z_{\theta}) = - \frac{1}{q+1} \, \e^{2
\frac{q-p}{p-1} -N}
\int_{\re^{N}} b\left(  \frac{x}{\e }  \right)  z_{0}^{q+1}(x+\theta
) dx,$$

\noindent and this expression of course vanishes as $\e \rightarrow
0$,
because the integral is bounded.
Hence, let us assume $N\geq 2 \frac{q-p}{p-1} $.  We obtain

$$ \e^{-N} G_{2}(\e , z_{\theta}) = - \frac{1}{q+1} \, \e^{2
\frac{q-p}{p-1} -N}
\int_{\re^{N}} b\left(  \frac{x}{\e }  \right)  z_{0}^{q+1}(x+\theta
) dx \leq $$

$$ C \e^{2 \frac{q-p}{p-1} -N} \, \left( \int_{\re^{N}} \Big\vert
b\left(  \frac{x}{\e }  \right)
\Big\vert^{\b } \, dx \right)^{\frac{1}{\b }} \,
\left( \int_{\re^{N}}  z_{0}^{\frac{\b (q+1)  }{\b -1}}  (x+
\theta )\, dx
\right)^{\frac{\b -1  }{\b }}\leq $$

$$ C \, \e^{2 \frac{q-p}{p-1} -N + \frac{N}{\b }} \left(
\int_{\re^{N}} \vert b(y)
\vert^{\b } \, dy \right)^{\frac{1}{\b }},$$

\noindent where $\b $ is given by ${\bf (b_{2})}$. This term goes to
zero since
$2 \frac{q-p}{p-1} -N + \frac{N}{\b } >0$. We have  now proved
(\ref{eq:lim4}), hence, by (\ref{eq:lim3}), (\ref{eq:lim1}) is also proved.

\medskip

\noindent Let us go to the proof of (\ref{eq:lim2}). Again we will
study separately $G_{1}'$ and $G_{2}'$.
We have

$$\vert \vert G_{1}'(\e , z_{\theta}) \vert \vert = \sup_{\vert \vert
v \vert \vert \leq 1} \big\vert
  (G_{1}'(\e , z_{\theta}) \vert v ) \big\vert =  \sup_{\vert
\vert v
\vert \vert \leq 1}
\Big\vert \int_{\re^{N}} \left( a\left( \frac{x}{\e} \right) -A
\right) \, z_{\theta}^{p}v dx \Big\vert \leq $$

$$\sup_{\vert \vert v \vert \vert \leq 1}  \left\{ \,  \left(
\int_{\re^{N}}
\Big\vert
a\left( \frac{x}{\e} \right) -A \Big\vert^{\frac{2N}{N+2}}  \,
z_{\theta}^{p\frac{2N}{N+2}} \, dx
\right)^{\frac{N+2}{2N}}   \, \left(  \int_{\re^{N}} \vert v
\vert^{2^{*}} dx \right)^{\frac{1}{2^{*}}} \right\}  \leq $$

$$C \e^{\frac{N}{2} + 1} \, \left( \int_{\re^{N}}  \vert
a( y) -A \vert^{\frac{2N}{N+2}} dy \right)^{\frac{N+2}{2N}} ,$$

\noindent hence

\begin{equation}\label{eq;lim5}
G_{1}'(\e , z_{\theta}) = O(\e^{\frac{N}{2} + 1} ).
\end{equation}

\noindent As to $G_{2}'(\e , z_{\theta})$ we obtain, arguing ad
before,

$$ \vert \vert G_{2}'(\e , z_{\theta})  \vert \vert \leq C
\e^{2\frac{q-p}{p-1} }
\left( \int_{\re^{N}} \Big\vert b\left( \frac{x}{\e} \right)
\Big\vert^{\frac{2N}{N+2}}
  z_{\theta}^{q \frac{2N}{N+2}} dx \right)^{\frac{N+2}{2N}} .$$

  \noindent By the usual change of variables $y=
\frac{x}{\e}$ and using ${\bf (b_{2})}$, we obtain

  $$G_{2}'(\e , z_{\theta}) = o(\e^{\frac{N}{2} + 1} ).$$

  \noindent From this and (\ref{eq;lim5}) we readily get
(\ref{eq:lim2}).

\end{pf}

\bigskip

\begin{Remark}
Notice that in the abstract results of section 2 the function
$\Gamma $ is defined on the manifold $Z$ of critical point of the
unperturbed functional $f_{0}$. In the present case this manifold is
diffeomorphic to $\re^{N}$, so we consider $\Gamma$ as a function
defined on $\re^{N}$.
\end{Remark}

\bigskip

\noindent We can now conclude the proof of theorem \ref{th:1}. We 
know that $z_{0}$
has a strict (global) maximum in $x= 0$, so that $\Gamma$ has a 
(strict) global maximum or minimum
(depending on the sign of $\int (a(y)-A)dy$) at $\theta =0$. We can 
then apply theorem 2.1,
setting $z^{*}= 0$ and, for example, $\delta =1$. We obtain a family
$\{ ( \e , u_{\e} )\}   \subset \re
\times
H^1 (\re^{N} ) $ such that $u_{\e}$ is a critical point of $f_{\e} $,
hence a solution of (2), and $\{  u_{\e} \}  $ is a bounded set in
$H^1 (\re^{N} ) $.
To be precise, we have

$$ u_{\e} (x)= z_{0} (x + \theta_{\e} ) + w(\e , \theta_{\e})(x) $$

\noindent where $\vert  \theta_{\e}\vert \leq 1 $ and $w(\e , \theta_{\e})
\rightarrow 0$ as $\e
\rightarrow 0 $. As $p< 1+ 4/N $, we obtain a family $(\l ,
\psi_{\l} )$
of solutions of (1) such that $\psi_{\l}  \rightarrow
0$ in $H^{1}(\re^{N} )$ as $\l \rightarrow 0$.

\bigskip

\begin{Remark}
The hypothesis ${  \bf (a_{2}) }$ is not used to prove the properties 
${\bf  (G_{0}-G_{3} )}$. It is used
to apply theorem 2.1, and in particular to say that there are $z^{*}, 
\delta$ such that (\ref{eq:nondeg}) holds.
If $\int (a(y)-A)dy = 0$ then $\Gamma$, as defined in 
(\ref{eq:gamdef}), is identically zero. It has
critical points, but of course they are not stable under perturbations,
so we can not conclude that they give rise to critical points of
$f_{\e}$.
\end{Remark}

\bigskip
\bigskip

\section{Second bifurcation result}
In this section we prove theorem 1.2. As before we have to prove that
$G, G' $ and $G''$ are continuous functions. Notice that in the proof
of ${ \bf (G_{0} )}$ and ${ \bf (G_{2} )}$ we will consider just the
function $G_{1}$, because the arguments of lemmas \ref{le:1} and
\ref{le:2} for the function
$G_{2}$ use only hypothesis ${
\bf (b_{1} )}$ which is unchanged. On the contrary in the proof of ${
\bf (G_{3} )}$ we will
study both $G_{1} $  and $G_{2}$. As above we assume $N\geq 3$ and
$1<p<q\leq \frac{N+2}{N-2}$. Let us see first that ${\bf (G_{0} )}$
is
satisfied.

\begin{Lemma}\label{le:3.1}
Assume ${\bf (a_{3} )}$. Then $G_{1}$ is
continuous.
\end{Lemma}
\begin{pf}
In the case  $(\e , u ) \rightarrow  (\e_{0} , u_{0} )$, $\e_{0}
\not= 0$ we can repeat word by word the arguments of lemma
\ref{le:1}, because there we
used only the fact that
$a$ is continuous and bounded, which is still true in the present
case. Hence, let us suppose
$(\e , u ) \rightarrow  (0 , u_{0} )$. Let us fix $\eta > 0$ and, by
${\bf (a_{3} ) }$, $M_{\eta } > 0$ such
that $\vert a(y)-A \vert < \eta $ if $\vert y \vert > M_{\eta } $. We
obtain

$$ (p+1) \vert G_{1} (\e , u ) \vert =  \Big\vert  \int_{\re^{N}
}\left( a \left( \frac{x}{\e} \right) -A \right)
\, \vert u \vert^{p+1} \, dx \Big\vert \leq $$

$$ \int_{\vert \frac{x}{\e } \vert \leq M_{\eta} } \Big\vert a \left(
\frac{x}{\e} \right) -A  \Big\vert
\, \vert u \vert^{p+1} dx + \int_{\vert \frac{x}{\e } \vert >
M_{\eta} } \Big\vert a \left( \frac{x}{\e} \right)
  -A  \Big\vert
\, \vert u \vert^{p+1} dx \leq $$

$$ C  \int_{\vert x  \vert \leq \e    M_{\eta} } \vert u \vert^{p+1}
dx + \eta \,  \int_{\re^{N} }
\vert u \vert^{p+1} dx \leq$$

$$C \int_{\re^{N}} \big\vert \vert u
\vert^{p+1} -\vert u_{0} \vert^{p+1} \big\vert dx + C
\int_{\vert x  \vert \leq \e    M_{\eta} } \vert u_{0} \vert^{p+1} dx
+\eta \,  \int_{\re^{N}}\vert u \vert^{p+1} dx
  .$$

\noindent As $\e \rightarrow 0$ and $ u \rightarrow u_{0}$ the first
two terms vanish, so

$$ \limsup_{ (\e , u ) \rightarrow  (0 , u_{0} ) } \vert G_{1} (\e ,
u ) \vert \leq C \eta .$$

\noindent This is true for any $ \eta > 0$, so we conclude $ G_{1}
(\e , u ) \rightarrow 0$ when $(\e , u )
  \rightarrow  (0 , u_{0} )$, and the lemma is proved.
\end{pf}

\bigskip
\noindent This proves that ${\bf (G_{1} )} $ holds. We want now to 
show that also ${\bf (G_{2} )} $ holds.

\bigskip

\begin{Lemma}
Assume ${\bf (a_{3} )}$. Then $G'_{1}$ and
$G''_{1}$ are continuous.
\end{Lemma}
\begin{pf}
We will show the continuity of $G'_{1}$, the argument for $G''_{1}$
is similar. Assume
$(\e , u ) \rightarrow  (\e_{0 }, u_{0} )$ with $\e_{0 } \not= 0$. In
this case one can argue
exactly as in lemma \ref{le:2} to obtain $\vert \vert G_{1}' (\e ,u )
-G_{1}' (\e_{0}, u_{0}) \vert
\vert \rightarrow 0$. Hence, let us now assume
$(\e , u ) \rightarrow  (0 , u_{0} )$. For each $\eta > 0$ let us fix
$M_{\eta } > 0$ as in the previous lemma. We obtain

$$ \vert \vert G'_{1}  (\e, u )  \vert \vert = \sup_{ \vert \vert v
\vert \vert \leq 1 } \Big\vert
\left( a \left( \frac{x}{\e} \right) -A \right) \, \vert u
\vert^{p-1} uv \, dx  \Big\vert \leq
$$

$$\sup_{ \vert \vert v \vert \vert \leq 1 } \left\{ \int_{ \vert
\frac{x}{\e} \vert  \leq M_{\eta} }
\Big\vert  a \left( \frac{x}{\e} \right) -A  \Big\vert  \, \vert u
\vert^{p} \vert v  \vert \, dx  \right\} +
\sup_{ \vert \vert v \vert \vert \leq 1 } \left\{ \int_{ \vert
\frac{x}{\e} \vert  >  M_{\eta} }\Big\vert
a \left( \frac{x}{\e} \right) -A \Big\vert \, \vert u \vert^{p} \vert
v \vert  \, dx \right\}  \leq
$$

$$C\sup_{ \vert \vert v \vert \vert \leq 1 } \left\{  \int_{ \vert   x
\vert \leq
\e  M_{\eta} } \vert u \vert^{p} \vert v
  \vert  \, dx \right\} + \eta
\sup_{ \vert \vert v \vert \vert \leq 1 } \left\{ \int_{ \re^{N} } \vert u
\vert^{p} \vert v \vert  \, dx  \right\}\leq
$$

$$
C\sup_{ \vert \vert v \vert \vert \leq 1 } \left\{  \left( 
\int_{\vert x \vert \leq \e M_{\eta}} \vert u
\vert^{p+1}dx \right)^{p/(p+1)} \, \left( \int_{\vert x \vert \leq \e 
M_{\eta}} \vert v \vert^{p+1} dx
\right)^{1/(p+1)} \right\}  + $$

$$\eta  \sup_{ \vert \vert v \vert \vert \leq 1 } \left\{  \left( 
\int_{\re^{N}} \vert u
\vert^{p+1}dx \right)^{p/(p+1)} \, \left( \int_{\re^{N}} \vert v \vert^{p+1} dx
\right)^{1/(p+1)} \right\} \leq
  $$

$$
C \left( \int_{\vert x \vert \leq \e M_{\eta}} \vert u
\vert^{p+1}dx \right)^{p/(p+1)} + C \eta \left(  \int_{\re^{N}} \vert u
\vert^{p+1}dx \right)^{p/(p+1)}.
$$

\noindent Arguing as before we then obtain

$$ \limsup_{(\e , u ) \rightarrow  (0 , u_{0} )}   \vert \vert
G'_{1}  (\e, u )  \vert \vert \leq C  \eta  $$

\noindent for all $\eta > 0$, hence $\lim_{(\e , u ) \rightarrow  (0
, u_{0} )}   \vert \vert G'_{1}  (\e, u )
  \vert \vert = 0$ and the lemma is proved.
\end{pf}

\bigskip

\noindent In the next lemma  we prove that ${ \bf (G_{3} ) } $ is
satisfied.

\bigskip

\begin{Lemma}
Assume ${\bf (a_{3} ) } $, ${\bf (b_{1} ) } $ and ${\bf (b_{3} ) } $.
Define

$$\Gamma (\theta )  = - \frac{L}{p+1}  \int_{\re^{N}} \vert x
\vert^{-\gamma } \, z_{0}^{p+1}
  ( x+\theta ) dx .
$$
\noindent Then, for all $\theta \in \re^{N}$, we have

\begin{equation}\label{eq:lim6}
\lim_{\e \rightarrow 0} \frac{G(\e , z_{\theta}) }{\e^{\gamma}}=
\Gamma (\theta )
\end{equation}

\noindent and

\begin{equation}\label{eq:lim7}
G'(\e , z_{\theta}) = o(\e^{\gamma / 2} ).
\end{equation}
\end{Lemma}
\begin{pf}
As usual we will study separately $G_{1}$ and $G_{2}$. We have

$$ G_{1}(\e , z_{\theta} )  = - \frac{1}{p+1}  \int_{\re^{N}} \left(
a\left( \frac{x}{\e}\right)   -A
   \right) z_{0}^{p+1} ( x+\theta ) dx =$$

$$ - \frac{\e^{\gamma}}{p+1}  \int_{\re^{N}} \left( a\left(
\frac{x}{\e}\right)   -A   \right)
\frac{\vert  x \vert^{\gamma}}{\e^{\gamma}}
\frac{ z_{0}^{p+1} ( x+\theta ) }{\vert  x \vert^{\gamma}} dx  .$$

\noindent We know that, for all $x\not= 0$,

$$ \left( a\left( \frac{x}{\e}\right)   -A   \right) \frac{\vert  x
\vert^{\gamma}}{\e^{\gamma}}
  \rightarrow L \quad \quad \hbox{as} \quad \quad \e \rightarrow 0,
$$

\noindent while, since $\gamma < N$ and $z_{0}$ has an exponential
decay at infinity,
$\vert  x \vert^{-\gamma}  \ z_{0} (x)  \in L^{1} (\re^{N}) .$ Hence
by dominated convergence we get

\begin{equation}\label{eq:lim8}
\lim_{\e \rightarrow 0} \frac{G_{1}(\e , z_{\theta}) }{\e^{\gamma}}=
\Gamma (\theta ) .
\end{equation}

\noindent To study $ \frac{1}{\e^{\gamma}} G_{2}(\e , z_{\theta})$ we 
can repeat the argument used in Lemma \ref{le:3},
distinguishing the cases
$\gamma < 2\frac{q-p}{p-1}$ and $\gamma \geq  2\frac{q-p}{p-1}$, and
using $\bf{(b_{3})}$ instead of ${\bf(b_{2})}$. We obtain

\begin{equation}\label{eq:lim9}
\lim_{\e \rightarrow 0} \frac{G_{2}(\e , z_{\theta}) }{\e^{\gamma}}= 0,
\end{equation}

\noindent  and (\ref{eq:lim6}) follows from (\ref{eq:lim8}) and
(\ref{eq:lim9}).
\par \noindent Let us now prove (\ref{eq:lim7}). We study first $
G'_{1}$ then $G'_{2}$. With the same
arguments of lemma \ref{le:3} we get

$$ \vert \vert G'_{1}(\e , z_{\theta})  \vert \vert  \leq C \left(
\int_{\re^{N}} \Big\vert a \left( \frac{x}{\e}
  \right) -A \Big\vert^{\frac{2N}{N+2}} z_{\theta}^{p \frac{2N}{N+2}
}(x) \, dx \right)^{ \frac{N+2}{2N}}
.$$

\noindent We have to distinguish three cases.
\bigskip

{\bf First case:}  $\gamma < \frac{N+2}{2}$.
\par  We obtain

$$ \left( \int_{\re^{N}} \Big\vert a \left( \frac{x}{\e} \right) -A
\Big\vert^{\frac{2N}{N+2}}
  z_{\theta}^{p \frac{2N}{N+2} }(x) \, dx \right)^{ \frac{N+2}{2N}}=$$

$$ \e^{\gamma} \left( \int_{\re^{N}} \Big\vert \left( a \left(
\frac{x}{\e} \right) -A \right)
\frac{\vert x \vert^{\gamma}}{\e^{\gamma} }
\Big\vert^{\frac{2N}{N+2}} \frac {z_{\theta}^{p \frac{2N}{N+2} }(x)}
{\vert x \vert^{\gamma \frac{2N}{N+2} }} \, dx \right)^{
\frac{N+2}{2N}}\leq $$

$$C \e^{\gamma} \left( \int_{\re^{N}}\frac{ z_{\theta}^{p
\frac{2N}{N+2} }(x) }
{  \vert x \vert^{\gamma \frac{2N}{N+2}}  } \, dx \right)^{
\frac{N+2}{2N}}\leq C \e^{\gamma},
$$

because $ \gamma \frac{2N}{N+2}<N$, hence $ z_{\theta}^{p
\frac{2N}{N+2} }(x) \, \vert x
\vert^{- \gamma \frac{2N}{N+2}} \in L^{1}(\re^{N}) $. Therefore in
this case

$$ \vert  \vert  G'_{1}(\e , z_{\theta})  \vert  \vert = O(
\e^{\gamma }) . $$

\bigskip

{\bf Second case:}  $\gamma > \frac{N+2}{2}$.
\par  In this case the function $\vert a(x) -A
\vert^{\frac{2N}{N+2}}$ is in $L^{1}(\re^{N})$, because it
is bounded and at infinity it is asymptotic to $\vert x \vert^{-
\gamma \frac{2N}{N+2}}$, and
$\gamma \frac{2N}{N+2}>N$. Therefore, by the usual change of
variables $y=\frac{x}{\e}$ we obtain

$$ \left( \int_{\re^{N}} \Big\vert a \left( \frac{x}{\e} \right)-A
\Big\vert^{\frac{2N}{N+2}} z_{\theta}^{p \frac{2N}{N+2} }(x) \, dx
\right)^{ \frac{N+2}{2N}}\leq $$

$$ C \e^{\frac{N+2}{2}} \left( \int_{\re^{N}}\vert a(y) -A \vert^{
\frac{2N}{N+2}}dy \right)^{ \frac{N+2}{2N}}\leq C \e^{\frac{N+2}{2}}.
$$
Hence, recalling that $\gamma <N$, we obtain
$$ \vert  \vert  G'_{1}(\e , z_{\theta})  \vert  \vert = O(
\e^{\frac{N}{2}+1 }) = o( \e^{\gamma / 2}). $$

\bigskip

{\bf Third case:}  $\gamma = \frac{N+2}{2}$.
\par In this case we apply H\"older inequality using as conjugate
exponents, instead of $\frac{2N}{N+2}$ and $\frac{2N}{N-2}$, any $s,
s'$
such that $s$ is smaller than  $\frac{2N}{N-2}$ but near to it, so
that $s'$ is bigger
than $\frac{2N}{N+2}$  but near to it. In this way we obtain

$$ \vert \vert   G'_{1}( \e , z_{\theta} )   \vert \vert =
\sup_{\vert \vert  v \vert \vert  \leq 1 }
\vert  (G'_{1}( \e , z_{\theta} ) \vert v ) \vert  =$$

$$\sup_{\vert \vert  v \vert \vert  \leq 1 } \Big\vert \int_{\re^{N}} \left( a
\left( \frac{x}{\e} \right)-A \right) z_{\theta}^{p} v dx \Big\vert \leq
\sup_{\vert \vert  v \vert \vert  \leq 1 } \left(  \int_{\re^{N}}
\Big\vert a \left( \frac{x}{\e} \right)-A
\Big\vert^{s'} z_{\theta}^{ps'}dx \right)^{1/s'} \, \left(
\int_{\re^{N}} \vert v \vert^{s}dx \right) ^{1/s}\leq $$

$$C \left(  \int_{\re^{N}}  \Big\vert a \left( \frac{x}{\e}
\right)-A
\Big\vert^{s'} dx \right)^{1/s'}.$$

\noindent It is $s' \gamma >N$, so, as before, $\vert a-A \vert^{s'}\in L^{1}
(\re^{N})$. We can then apply the usual change of variables to obtain

$$ \vert \vert   G'_{1}( \e , z_{\theta} )   \vert \vert \leq C
\e^{N/s'}.$$

We have that $s'$ is near $\frac{2N}{N+2}= \frac{N}{\gamma}$, so
$N/s'$ is near $\gamma$. Hence we obtain

$$ \vert \vert   G'_{1}( \e , z_{\theta} )   \vert \vert  = O(
\e^{N/s'})= o(\e^{\frac{\gamma}{2}}).$$

\bigskip

\noindent We have concluded the study of $ G'_{1}( \e , z_{\theta}
)$. As to $G'_{2}( \e , z_{\theta})$, the same argument of lemma 
\ref{le:3} gives

$$ \vert \vert   G'_{2}( \e , z_{\theta} )   \vert \vert  =
o(\e^{\frac{N}{2}+1})=
o(\e^{\frac{\gamma}{2}}).$$

\noindent In  this way the lemma
is completely proved.
\end{pf}

\bigskip

\noindent We want now to complete the proof of Theorem 1.2. As in the
previous section, we have only to prove that the function $\Gamma$
satisfies the hypotheses of theorem 2.1. Let us prove that

\begin{equation}\label{eq:Gam1}
\hbox{there is}\,\, R>0 \,\, \hbox{ such that either} \,\,
\min_{\vert \theta \vert = R} \Gamma (\theta ) > \Gamma (0)
\,\, \hbox{ or} \,\,
\max_{\vert \theta \vert = R} \Gamma (\theta ) < \Gamma (0).
\end{equation}

\noindent To prove (\ref{eq:Gam1}) we first notice that $\Gamma$ is 
continuous and $\Gamma (\theta )$ is either
positive on all $\re^{N}$ or negative on all $\re^{N}$. Then we claim that

\begin{equation}\label{eq:Gam2}
\lim_{\vert  \theta \vert  \rightarrow  +\infty } \Gamma (\theta ) =
0.
\end{equation}

\noindent To prove (\ref{eq:Gam2}), let us write

$$\vert   \Gamma (\theta ) \vert = C \int_{\re^{N}} \vert x
\vert^{-\gamma } \, z_{0}^{p+1}
  ( x+\theta ) dx=
C \int_{\vert x \vert \leq 1} \vert x \vert^{-\gamma } \, z_{0}^{p+1}
  ( x+\theta ) dx + C \int_{\vert x \vert > 1} \vert x \vert^{-\gamma
}
\, z_{0}^{p+1}
  ( x+\theta ) dx .$$

\noindent It is $\vert x \vert^{-\gamma } \in L^{1} (B_{1}) $, while
$ z_{0}^{p+1}
  ( x+\theta )\rightarrow 0$ as $\vert  \theta \vert \rightarrow
+\infty $, for all $x$, so by dominated convergence the first
integral vanishes as $\vert  \theta \vert \rightarrow +\infty$. For
the
second integral we write

$$\int_{\vert x \vert > 1} \vert x \vert^{-\gamma } \, z_{0}^{p+1}
  ( x+\theta ) dx = \int_{\vert y-\theta  \vert > 1} \vert y-\theta
\vert^{-\gamma } \, z_{0}^{p+1}
  ( y) dy = \int_{\re^{N} } \chi_{\theta}(y) \vert y-\theta
\vert^{-\gamma } \, z_{0}^{p+1}
  ( y) dy, $$

\noindent where $\chi_{\theta}$ is the characteristic function of
the set $\{  y\in \re^{N} \, \vert  \,
\vert y-\theta  \vert > 1 \}$. It is trivial to see that
$$\chi_{\theta}(y) \vert y-\theta \vert^{-\gamma } \leq 1$$
\noindent for all $y, \theta \in \re^{N}$ ($y\not= \theta$), and that
$\chi_{\theta}(y) \vert y-\theta \vert^{-\gamma } \rightarrow 0$ as
$\vert \theta \vert \rightarrow +\infty$. Again by dominated
convergence we obtain that also the second integral vanishes when
$\vert \theta \vert \rightarrow +\infty$. So (\ref{eq:Gam2}) is
proved, hence also
(\ref{eq:Gam1}). We can apply theorem 2.1 and argue as in the
previous section.

\section{Continuous branches of solutions.}
In this section we prove that in some cases the families of solutions
bifurcating from $(0,0)$, that we have found in the previous
sections, form a curve. We first will prove some abstract results
(following the frame of section 2), then we will apply these result
to problems (\ref{eq:gene}) and  (\ref{eq:gene1}). So let us come
back to the abstract frame of section 2. To make easier the passage
from the abstract frame to the applications, we will denote by
$z_{\theta}$, $\theta \in \re^{d}$, the elements of $Z$. Notice that
our arguments are local in nature and we will work in the neighborhood
of a fixed point, so we can assume, without loss of generality,
that the manifold $Z$ is given by a unique map $\theta \rightarrow
z_{\theta}$. We will indicate with $\partial_{i}z$, $\partial_{ij}z$
the derivatives of $z_{\theta}$ with respect to the parameter
$\theta$, that is
$$\partial_{i}z = \frac{\partial z}{\partial \theta_{i}}(\theta),
\quad \quad \partial_{ij}z = \frac{\partial^{2} z}{\partial \theta_{i}
\partial \theta_{j} }(\theta)  . $$

\noindent About the manifold $Z$ of critical points we will also
assume the following hypothesis, which is satisfied in our applications.

\bigskip

\begin{description}
        \item[${\bf (H)}$] $ (\partial_{i}z \, \vert \, 
\partial_{j}z) =0$ if $i\not= j$,
        $\vert \vert  \partial_{i}z  \vert \vert = c$ (independent of 
$i$ and $\theta$),
        $ (\partial_{ij}z \, \vert \, \partial_{l}z) =0$ for all 
$i,j,l =1,..,d$.
\end{description}

\bigskip

\noindent About the functionals $F,G$ we will assume two different
types of hypotheses. Recall that $\a, \Gamma$ are those given in 
hypothesis ${\bf (G_3)}$.

\begin{description}
        \item[${\bf (F_4)}$] $F$ is of class $C^{4}$;

        \item[${\bf (G_4)}$]  $G$ is of class $C^{4}$ with respect to 
$u$ and the map $(\e ,
u) \rightarrow G''' (\e , u)$ is continuous;

        \item[${\bf (G_5)}$]
$\Gamma$ is $C^{2}$ and, if $\theta_{\e}$ is a family such that 
$\theta_{\e} \rightarrow \theta$ as
$\e \rightarrow 0$,
then
$$  \lim_{\e \rightarrow 0} \left[  \e^{-\a} (G'(\e , z_{\theta_{\e}} ) \,
\vert  \, \partial_{ij} z_{\theta_{\e}} ) +\e^{-\a} (G''(\e ,
z_{\theta_{\e}} )\partial_{i} z_{\theta_{\e}} \, \vert  \,
\partial_{j} z_{\theta_{\e}} ) \right] =
\partial_{ij}\Gamma (\theta ) .$$

$$
  \lim_{\e \rightarrow 0} \e^{-\a /2} G'' (\e ,
z_{\theta_{\e}})\partial_{i} z_{\theta_{\e}} =0 ,
\qquad
\lim_{\e \rightarrow 0} \e^{-\a /2} G'' (\e ,
z_{\theta_{\e}})\partial_{ij} z_{\theta_{\e}} =0 $$

$$
\textrm{and}\qquad
\lim_{\e \rightarrow 0} \e^{-\a /2} G''' (\e ,
z_{\theta_{\e}})[ \partial_{i} z_{\theta_{\e}} ,\partial_{j}
z_{\theta_{\e}}] =0 .$$

\end{description}

\bigskip

\begin{description}
        \item[${\bf (F_4)'}$] $F$ is of class $C^{3}$;

        \item[${\bf (G_4)'}$]$G$ is of class $C^{3}$ with respect to 
$u$ and the map $(\e ,
u) \rightarrow G''' (\e , u)$ is continuous.

        \item[${\bf (G_5)'}$]  $G'(\e , u ) = O(\e^{\a  }) $ for all 
$u\in E$, $\Gamma$ is $C^{2}$ and,
        if $\theta_{\e}$ is a family such that $\theta_{\e} 
\rightarrow \theta$ as $\e \rightarrow 0$,
then
$$  \lim_{\e \rightarrow 0} \left[ \e^{-\a} (G'(\e , z_{\theta_{\e}} ) \,
\vert  \, \partial_{ij} z_{\theta_{\e}} ) +\e^{-\a} (G''(\e ,
z_{\theta_{\e}} )\partial_{i} z_{\theta_{\e}} \, \vert  \,
\partial_{j} z_{\theta_{\e}} ) \right] =
\partial_{ij}\Gamma (\theta ) ,$$

$$
\lim_{\e \rightarrow 0} \e^{-\a /2} G'' (\e ,
z_{\theta_{\e}})\partial_{i} z_{\theta_{\e}} =0 .
$$

\end{description}

\bigskip

\noindent We can now prove two abstract theorems.
\bigskip

\begin{Theorem}\label{th:abstract1}
Assume ${ \bf (H), (F_{0}-F_{4}) }$ and ${ \bf (G_{0}-G_{5})}$. For a 
given $\theta
\in \re^{d}$, and for any small
  $\e$'s, let us suppose that there is a critical point
$ u_{\e } \in Z_{\e}$ of
$f_{\e}$,
such that $u_{\e }= z_{\theta_{\e}} + w(\e , \theta_{\e})$ and
$\theta_{\e} \rightarrow \theta$ as
  $\e \rightarrow 0$. Assume that $z_{\theta}= \lim_{\e}
z_{\theta_{\e}} $ is nondegenerate for
  the restriction of $f_{0}$ to $(T_{z_{\theta}}Z)^{\perp}$, with
Morse index equal to $m_{0}$, and that
  the hessian matrix $D^{2}\Gamma (\theta )$ is positive or negative
definite.
\par Then $u_{\e }$, for small $\e $'s , is a nondegenerate critical 
point for $f_{\e }$
with Morse index equal to $m_{0}$
if $D^{2}\Gamma (\theta)$ is positive, to $m_{0} + d$ if $D^{2}\Gamma
(\theta )$ is negative. A
a consequence, the critical points of  $f_{\e }$ form a continuous
curve.
\end{Theorem}
\begin{pf}
Let us write
\begin{equation}\label{eq:dec}
E= E^{+}\oplus E^{0} \oplus E^{-}
\end{equation}

\noindent where $E^{0}= T_{z_{\theta}}Z$, dim$(E^{-})=m_{0}$ and there
exists $\delta >0$ such that

$$
\left\{
          \begin{array}{rcl}
D^2f_0(z_{\theta})[v,v] & >\delta \vert \vert  v \vert \vert^{2} & 
\forall\; v\in E^{+},\\
D^2f_0(z_{\theta})[v,v] & <-\delta \vert \vert v \vert \vert^{2} 
&\forall\; v\in E^{-}.
            \end{array}
\right.
$$

\noindent From the hypothesis $f_{0}'(z_{\eta } )= 0$ for all $\eta
\in \re^{N}$ it is easy to deduce

$$ D^{2} f_{0} (z_{\theta }) [\partial_{i} z_{\theta} , \partial_{j}
z_{\theta}] =0 .$$

\noindent Let us define
$$\varphi_{i}^{0}  = \frac{1}{\vert \vert     \partial_{i}
z_{\theta}     \vert \vert } \,  \partial_{i} z_{\theta}.$$

\noindent The set $\{  \varphi_{i}^{0} \}_{i=1,..,d}$ is an
orthonormal base for $E^{0}$. Let $\l_{1},\ldots , \l_{d}$ be the 
eigenvalues of
the symmetric matrix $D^{2} f_{0} (z_{\theta })$ on $E^{-}$. Of 
course $\l_{i}<0$ for all $i$, and let
$\l_{0}= \max_{i}\l_{i}<0$. Let $\{ t_{i}^{0} \}_{i=1,..,m_{0}}$ be
an orthonormal base for $E^{-}$ such that $D^{2}f_{0} [t_{i}^{0} , 
t_{j}^{0} ] = 0$ if $i\not= j$,
$D^{2}f_{0} [t_{i}^{0} , t_{i}^{0} ] =\l_{i}$. By orthogonality of
the decomposition (\ref{eq:dec}),
we have $(  \varphi_{i}^{0} \, \vert \, t_{j}^{0} ) =0$ for all
$i,j$. Define

$$ \varphi_{i}^{\e} =  \frac{1}{\vert \vert     \partial_{i}
z_{\theta_{\e}}     \vert \vert } \,  \partial_{i} z_{\theta_{\e}}.$$

\noindent The set $\{ \varphi_{i}^{\e} \}_{i=1,..,d }$ is an 
orthonormal base for the tangent
space $T_{z_{\theta_{\e}}}Z$,
space that we denote $E^{0}_{\e }$. Notice that $ \varphi_{i}^{\e}
\rightarrow  \varphi_{i}^{0}$ as $\e \rightarrow 0$.
\par \noindent For $i= 1,..,m_{0}$ we want to find $\tau_{i}^{\e }$
such that, setting $t_{i}^{\e } =t_{i}^{0} + \tau_{i}^{\e }$, we
obtain, for all $i,j$,

\begin{equation}\label{eq:perp}
( t_{i}^{\e } \, \vert  \,  \varphi_{j}^{\e} )=0
\end{equation}

That is, we want

$$ 0 =( t_{i}^{\e } \, \vert  \,  \varphi_{j}^{\e} )= (t_{i}^{0} +
\tau_{i}^{\e } \, \vert  \,
\varphi_{j}^{0}+ (  \varphi_{j}^{\e}-  \varphi_{j}^{0})) =
(t_{i}^{0}  \, \vert \,  \varphi_{j}^{\e}-  \varphi_{j}^{0}) +
(\tau_{i}^{\e } \, \vert \,  \varphi_{j}^{\e} ), $$

\noindent hence

$$  (\tau_{i}^{\e } \, \vert \,  \varphi_{j}^{\e} )= - (t_{i}^{0}  \,
\vert \,  \varphi_{j}^{\e}-  \varphi_{j}^{0}).$$

\noindent So, we define
$$ \tau_{i}^{\e }= \sum_{j=1}^{d} - (t_{i}^{0}  \, \vert \,
\varphi_{j}^{\e}-  \varphi_{j}^{0})\varphi_{j}^{\e}  ,$$

\noindent and (\ref{eq:perp}) holds. Notice that $\tau_{i}^{\e}
\rightarrow 0$ as $\e \rightarrow 0$, so that
$$t_{i}^{\e } \rightarrow t_{i}^{0} $$

\noindent as $\e \rightarrow 0$. As $ \{ t_{i}^{0} \}_{i}$ is an
orthonormal base, the vectors $\{ t_{i}^{\e } \}$ are
linearly independent, for small $\e$'s. Let us define $E_{\e }^{-}$
the $m_{0}$-dimensional space spanned by $\{ t_{i}^{\e } \}$. For $v
\in E_{\e }^{-}$, $\vert \vert  v \vert \vert =1$, we have
$v= \sum_{k=1}^{m_{0}} \beta_{k} t_{k}^{\e }$ and we can write

$$ D^{2} f_{\e} (u_{\e }) [ v , v ] =
\sum_{l,k=1}^{m_{0}}\beta_{l} \beta_{k}  D^{2} f_{\e} (u_{\e }) [
t_{l}^{0 } +\tau_{l}^{\e }  , t_{k}^{0 } +\tau_{k}^{\e }] =
\sum_{l,k=1}^{m_{0}}\beta_{l} \beta_{k}  D^{2} f_{\e} (u_{\e }) [
t_{l}^{0 }   , t_{k}^{0 } ]+ o(1),  $$

\noindent where $o(1)$ vanishes as $\e \rightarrow 0$, uniformly in
$v$. By hypotheses ${\bf (F_{0} )}$, ${\bf (G_{1} )}$, ${\bf (G_{2}
)}$, we obtain
$$ D^{2} f_{\e} (u_{\e }) [ t_{l}^{0 }   , t_{k}^{0 } ] \rightarrow
D^{2} f_{0} (z_{\theta}) [ t_{l}^{0 }   , t_{k}^{0 } ].$$

\noindent As $t_{k}^{\e} \rightarrow t_{k}^{0}$, for $\e \rightarrow 
0$, and $ \{ t_{k}^{0} \} $ is orthonormal,
it is easy to see that, for small $\e$'s, $\vert \vert v \vert \vert 
=1$ implies $\sum_{k=1}^{m_{0}} \b_{k}^{2}
\geq \frac{1}{2}$. Hence

$$ D^{2} f_{\e} (u_{\e }) [ v , v ]= \sum_{k=1}^{m_{0}} \l_{k} 
\b_{k}^{2}  + o(1) \leq  \frac{\l_{0}}{2}
+o(1), $$

\noindent where $o(1) \rightarrow 0$ as $\e \rightarrow 0$, uniformly in $v$ if
$\vert \vert v \vert \vert =1$. Hence, for small $\e $,  $ D^{2}
f_{\e} (u_{\e })$ is negative definite in $E_{\e }^{-}$.
\bigskip
\noindent We now define

$$ E_{\e }^{+}= (E_{\e }^{0} \oplus E_{\e }^{-})^{\perp} ,$$

\noindent so that

$$E = E_{\e }^{+}\oplus E_{\e }^{0}\oplus E_{\e }^{-}.$$

\noindent We want now to prove that $ D^{2} f_{\e} (u_{\e })$ is
positive definite on $ E_{\e }^{+}$, for small $\e$. Let $P^{+}$ be
the orthogonal projection of $E$ to $E^{+}$.
\par \noindent We claim that there are $\delta_{0} >0$ and $\e_{0} >0$ such
that for all $\vert  \e \vert <
\e_{0}$ and all $v \in E_{\e }^{+}$, $\vert \vert v  \vert \vert 
=1$, it holds

$$ D^{2} f_{\e} (u_{\e })[v , v] > \delta_{0} .$$

\noindent We argue by  contradiction. If the claim is not true, then
there are sequences $\{ \e_{k} \} $, $\{ v_{k} \} \subset E_{\e_{k}
}^{+}$, with $\vert \vert v_{k} \vert \vert =1$ and $\e_{k} \rightarrow 0$ as
$k\rightarrow \infty$, such that

$$\frac{1}{k}> D^{2} f_{\e_{k} } (u_{\e_{k}  })[v_{k}   , v_{k}] =$$
  \begin{equation}\label{eq:abs1}
D^{2} f_{\e_{k} } (u_{\e_{k}  })[ P^{+}v_{k}   , P^{+}v_{k}] + 2D^{2}
f_{\e_{k} } (u_{\e_{k}  })[ P^{+}v_{k}   , v_{k}  - P^{+}v_{k}]+
D^{2} f_{\e_{k} } (u_{\e_{k}  })[v_{k}  - P^{+}v_{k}   , v_{k}  -
P^{+}v_{k}].
\end{equation}

\noindent  We recall that

$$ v_{k}  - P^{+}v_{k}= \sum_{i=1}^{m_{0}} (v_{k} \, \vert  \,
t_{i}^{0}) t_{i}^{0}+ \sum_{i=1}^{d}(v_{k} \, \vert \,
\varphi_{i}^{0} ) \varphi_{i}^{0},$$

\noindent and that, since $v_{k} \in E_{\e_{k}}^{+} $, $( v_{k} \, \vert
\,  t_{i}^{\e_{k}} ) =0$ and
$ (v_{k} \, \vert \, \varphi_{i}^{\e_{k}})=0$. Hence we have

$$ ( v_{k} \, \vert  \,  t_{i}^{0} )= (v_{k} \, \vert  \, t_{i}^{0} -
t_{i}^{\e_{k}} )\rightarrow 0 $$

\noindent as $k\rightarrow \infty$, because $ 
t_{i}^{\e_{k}}\rightarrow t_{i}^{0}$ and
$\{ v_{k} \}$ is bounded. In the same way we get

$$ ( v_{k} \, \vert  \,  \varphi _{i}^{0} ) \rightarrow 0 , $$

\noindent hence

\begin{equation}\label{eq:pro1}
  v_{k}  - P^{+}v_{k} \rightarrow 0 .
\end{equation}

\noindent Now (\ref{eq:abs1}) becomes

$$ \frac{1}{k}> D^{2} f_{\e_{k} } (u_{\e_{k}  })[ P^{+}v_{k}   ,
P^{+}v_{k}] + o(1) =$$
\begin{equation}\label{eq:abs2}
D^{2} f_{0 } (z_{\theta  })[ P^{+}v_{k}   , P^{+}v_{k}] +
(D^{2} f_{\e_{k} } (u_{\e_{k}  }) - D^{2} f_{0 } (z_{\theta  })) [
P^{+}v_{k}   , P^{+}v_{k}] + o(1).
\end{equation}

\noindent Thanks to the continuity hypotheses we have
$$(D^{2} f_{\e_{k} } (u_{\e_{k}  }) - D^{2} f_{0 } (z_{\theta  })) [
P^{+}v_{k}   , P^{+}v_{k}]= o(1),$$

\noindent while

$$ D^{2} f_{0} (z_{\theta  })  [ P^{+}v_{k}   , P^{+}v_{k}] > \delta
\vert \vert P^{+}v_{k} \vert \vert^{2}  ,$$

\noindent because $ P^{+}v_{k} \in E^{+}$. By (\ref{eq:pro1}) we
also obtain

$$ \vert \vert P^{+}v_{k} \vert  \vert \rightarrow  1 .$$

\noindent Hence (\ref{eq:abs2}) gives

$$ \frac{1}{k}> \delta + o(1), $$
\noindent a contradiction. So the claim is proved.
\par \noindent Up to now we have shown that, for small $\e$,
$D^{2}f_{\e } (u_{\e })$ is negative definite on $E^{-}_{\e} $ and
positive definite on $E^{+}_{\e} $. We want now to study the behavior
of $D^{2}f_{\e } (u_{\e })$ on $E^{0}_{\e} $. We will prove that
$D^{2}f_{\e } (u_{\e })$ is positive or negative definite accordingly with
$D^{2}\Gamma (\theta )$, and this will conclude the proof.
\par \noindent As first thing we recall that we have

$$
D^{2}f_{\e } (u_{\e } )[ \partial_{i} z_{\theta_{\e}}, \partial_{j}
z_{\theta_{\e}}] =
(\partial_{i} z_{\theta_{\e}}    \, \vert \,   \partial_{j}
z_{\theta_{\e}}) -
( F'' (u_{\e})  \partial_{i} z_{\theta_{\e}} \, \vert \,  \partial_{j}
z_{\theta_{\e}})  +
(G'' (\e , u_{\e})  \partial_{i} z_{\theta_{\e}} \, \vert \, \partial_{j}
z_{\theta_{\e}}).
$$

\noindent As $\partial_{i} z_{\theta_{\e}} \in $ ker$[I_{E} -
F''(z_{\theta_{\e}}) ]$ and
$w(0, z_{\theta_{\e}})=0$, developing $F''(u_{\e})$ and $G''(\e ,
u_{\e})$ and setting $w_{\e}= w(\e , \theta_{\e})$, we obtain

$$ D^{2}f_{\e } (u_{\e }) [ \partial_{i} z_{\theta_{\e}},
\partial_{j}  z_{\theta_{\e}}] =$$

$$ (\partial_{i} z_{\theta_{\e}}    \, \vert \,   \partial_{j}
z_{\theta_{\e}}) -
( F'' (z_{\theta_{\e}})  \partial_{i} z_{\theta_{\e}} \, \vert \, 
\partial_{j}  
z_{\theta_{\e}})  -
(F''' (z_{\theta_{\e}}) [ \partial_{i} z_{\theta_{\e}}, \partial_{j}
z_{\theta_{\e}}] \, \vert \,   w_{\e } )+$$
$$(G'' (\e , z_{\theta_{\e}})  \partial_{i} z_{\theta_{\e}} \, \vert 
\, \partial_{j}
z_{\theta_{\e}}) +
(G''' (\e , z_{\theta_{\e}}) [ \partial_{i} z_{\theta_{\e}},
\partial_{j}  z_{\theta_{\e}}] \, \vert \,   w_{\e } )+ O(\vert
\vert  w_{\e } \vert \vert^{2} ) =$$

$$ -(F''' (z_{\theta_{\e}}) [ \partial_{i} z_{\theta_{\e}},
\partial_{j}  z_{\theta_{\e}}] \, \vert \,   w_{\e } )+
(G'' (\e , z_{\theta_{\e}})  \partial_{i} z_{\theta_{\e}} \, \vert \, 
\partial_{j}
z_{\theta_{\e}}) +$$

$$ (G''' (\e , z_{\theta_{\e}}) [ \partial_{i} z_{\theta_{\e}},
\partial_{j}  z_{\theta_{\e}}] \, \vert \,   w_{\e } )+ O(\vert
\vert  w_{\e } \vert \vert^{2} ).$$

\noindent We have $ w_{\e } = o( \e^{\a / 2}) $ (see theorem 2.1),
hence from $\bf{ (G_{5}) }$ we deduce
$$ (G''' (\e , z_{\theta_{\e}}) [ \partial_{i} z_{\theta_{\e}},
\partial_{j}  z_{\theta_{\e}}] \, \vert \,   w_{\e } )= o(\e^{\a})$$

\noindent so that

\begin{equation}\label{eq:dev2}
D^{2}f_{\e } (u_{\e } )[ \partial_{i} z_{\theta_{\e}}, \partial_{j}
z_{\theta_{\e}}] =
-(F''' (z_{\theta_{\e}}) [ \partial_{i} z_{\theta_{\e}},
\partial_{j}  z_{\theta_{\e}}] \, \vert \,   w_{\e } )+
(G'' (\e , u_{\e})  \partial_{i} z_{\theta_{\e}} \, \vert \, \partial_{j}
z_{\theta_{\e}}) + o(\e^{\a}).
\end{equation}

\noindent By (\ref{eq:w}) we have

$$ z_{\theta_{\e }} + w_{\e} - F' (z_{\theta_{\e }} + w_{\e}) +G'(\e
, z_{\theta_{\e }} + w_{\e})=
\sum_{l} a_{l} \partial_{l} z_{\theta_{\e }}.$$

\noindent Developing $F'$ and $G'$ we obtain

$$ z_{\theta_{\e }} + w_{\e} - F' (z_{\theta_{\e }} )  - F''
(z_{\theta_{\e }} ) w_{\e} + G'(\e , z_{\theta_{\e }}) +G''(\e ,
z_{\theta_{\e }}) w_{\e} + O(\vert \vert w_{\e } \vert \vert^{2} )=
\sum_{l} a_{l} \partial_{l} z_{\theta_{\e }}.$$

\noindent By a scalar product with $\partial_{ij} z_{\theta_{\e}}$,
recalling ${\bf (H) }$ and the fact that $ z_{\theta_{\e}}=
F'(z_{\theta_{\e}})$, we get
$$(w_{\e} \,  \vert   \,  \partial_{ij} z_{\theta_{\e}}) - (F''
(z_{\theta_{\e}})w_{\e} \,  \vert   \, \partial_{ij} z_{\theta_{\e}})
+ (G' (\e , z_{\theta_{\e}}) \,  \vert   \, \partial_{ij}
z_{\theta_{\e}})+$$

$$(G'' (\e , z_{\theta_{\e}})w_{\e} \,  \vert   \, \partial_{ij}
z_{\theta_{\e}}) + O( \vert \vert  w_{\e }  \vert \vert^{2} ) =0.$$

\noindent By ${\bf (G_{5} ) }$ it is

$$ (G'' (\e , z_{\theta_{\e}})w_{\e} \,  \vert   \, \partial_{ij}
z_{\theta_{\e}}) = (G'' (\e , z_{\theta_{\e}}) \partial_{ij}
z_{\theta_{\e}} \,  \vert   \, w_{\e}) = o(\e^{\a}),$$

\noindent so we obtain

\begin{equation}\label{eq:dev1}
(w_{\e} \,  \vert   \,  \partial_{ij} z_{\theta_{\e}}) - (F''
(z_{\theta_{\e}})w_{\e} \,  \vert   \, \partial_{ij} z_{\theta_{\e}})
+ (G' (\e , z_{\theta_{\e}}) \,  \vert   \, \partial_{ij}
z_{\theta_{\e}})+
o(\e^{\a } ) = 0  .
\end{equation}

\noindent Deriving twice the equation $ z_{\eta} = F'(z_{\eta })$ with
respect to $\eta \in \re^{N}$ and computing the result at $\eta=
\theta_{\e}$ we obtain

$$ \partial_{ij} z_{\theta_{\e}} = F''(z_{\theta_{\e}})
\partial_{ij} z_{\theta_{\e}} +
F''' (z_{\theta_{\e}}) [  \partial_{i} z_{\theta_{\e}},  \partial_{j}
z_{\theta_{\e}} ] . $$

\noindent By a scalar product with $w_{\e }$ we get

$$ (  \partial_{ij} z_{\theta_{\e}} \,  \vert   \, w_{\e }) = (F''
(z_{\theta_{\e}}) \partial_{ij} z_{\theta_{\e}}) \,  \vert   \, w_{\e
}) +
(F''' (z_{\theta_{\e}}) [  \partial_{i} z_{\theta_{\e}},
\partial_{j} z_{\theta_{\e}} ] \,  \vert   \, w_{\e }). $$

\noindent Substituting this last identity in (\ref{eq:dev1}) we obtain

\begin{equation}\label{eq:dev3}
  - (F''' (z_{\theta_{\e}}) [  \partial_{i} z_{\theta_{\e}},
\partial_{j} z_{\theta_{\e}} ] \,  \vert   \, w_{\e })= ( G' (\e ,
z_{\theta_{\e}}) \,   \vert    \,   \partial_{ij} z_{\theta_{\e}}  )
+ o(\e^{\a}).
\end{equation}

\noindent From this and (\ref{eq:dev2}) we have

$$D^{2}f_{\e } (u_{\e } ) [ \partial_{i} z_{\theta_{\e}}, \partial_{j}
z_{\theta_{\e}}] =
( G' (\e , z_{\theta_{\e}}) \,   \vert    \,   \partial_{ij}
z_{\theta_{\e}}  ) + (G'' (\e , u_{\e})  \partial_{i} z_{\theta_{\e}} 
\, \vert \, \partial_{j}
z_{\theta_{\e}}) + o(\e^{\a}).$$

\noindent Hence, dividing by $\e^{\a}$, passing to the limit and
using hypothesis ${\bf (G_5)}$ we obtain

$$\lim_{\e \rightarrow 0} \frac{1}{\e^{\a}} D^{2}f_{\e } (u_{\e } )[
\partial_{i} z_{\theta_{\e}}, \partial_{j}  z_{\theta_{\e}}] = \
\lim_{\e \rightarrow 0}\frac{1}{\e^{\a}}    ( G' (\e ,
z_{\theta_{\e}}) \,   \vert    \,   \partial_{ij} z_{\theta_{\e}}  )
+ \frac{1}{\e^{\a}}  (G'' (\e , u_{\e})  \partial_{i} z_{\theta_{\e}} 
\, \vert \, \partial_{j}
z_{\theta_{\e}}) = \partial_{ij}
\Gamma (\theta ).$$

\noindent As $D^{2}\Gamma (\theta )$ is definite, so is $D^{2}f_{\e }
(u_{\e } )[ \partial_{i} z_{\theta_{\e}}, \partial_{j}
z_{\theta_{\e}}]$, for small $\e$'s.

\par \noindent Let us recall what we have proved up to now. Let us 
assume that $D^{2}\Gamma (\theta )$ is
definite positive (the other case is analogous). We have proved that 
there is a constant $ \delta > 0$ such that

$$ D^{2}f_{\e } (u_{\e } ) [v^{-}, v^{-} ] \leq - \delta \vert \vert 
v^{-} \vert \vert^{2}
\quad \quad \hbox{for all} \quad v^{-} \in E^{-}_{\e}, $$

$$ D^{2}f_{\e } (u_{\e } ) [v^{0}, v^{0} ] \geq  \delta \e^{\a} \vert 
\vert v^{0} \vert \vert^{2}
\quad \quad \hbox{for all}\quad v^{0} \in E^{0}_{\e}, $$

$$ D^{2}f_{\e } (u_{\e } ) [v^{+}, v^{+} ] \geq  \delta \vert \vert 
v^{+} \vert \vert^{2}
\quad \quad \hbox{for all}\quad v^{+} \in E^{+}_{\e}. $$

\noindent To conclude the proof we have to show that $ D^{2}f_{\e } 
(u_{\e } )$ is positive definite in
$E^{+}_{\e} + E^{o}_{\e}$. This does not derive directly from the 
previous statements because it is not
true, in general, that $ D^{2}f_{\e } (u_{\e } ) [v^{+}, v^{0}] = 0$. 
However, thanks to $ \bf{(G_{5})}$,
we have, for any $v^{+} \in E^{+}_{\e}$, $v^{0} \in E^{0}_{\e}$,
$$ \big\vert D^{2}f_{\e } (u_{\e } ) [v^{+}, v^{0}] \big\vert \leq 
o(\e^{\a/2} ) \vert \vert v^{+} \vert \vert \,
\vert \vert v^{0} \vert \vert .$$

\noindent Hence, for small $\e$'s and a suitable $\delta_{1}>0$, we obtain

$$ D^{2}f_{\e } (u_{\e } ) [v^{+}+ v^{0}, v^{+}+ v^{0}] \geq \delta 
\vert \vert v^{+} \vert \vert +
\delta  \e^{\a}\vert \vert v^{0} \vert \vert - o(\e^{\a/2})\vert 
\vert v^{+} \vert \vert \,
\vert \vert v^{0} \vert \vert  \geq \delta_{1} \vert \vert v^{+}+ 
v^{0} \vert \vert^{2} .$$

\noindent The proof is now complete.

\end{pf}

\bigskip

\noindent With small changes in the previous arguments one can prove
the following  theorem.

\begin{Theorem}\label{th:abstract2}
Assume ${   \bf   (H), (F_{0}-F_{3})   } $, ${   \bf (G_{0}-G_{3}) 
}$, ${   \bf (F_{4})' } $,
${   \bf (G_{4})'  }$,  ${   \bf (G_{5})'  }$. For a given $\theta 
\in \re^{d}$, and for
any small $\e$'s, let us suppose that there is a critical point
$u_{\e } \in Z_{\e}$ of
$f_{\e}$, such that $u_{\e }= z_{\theta_{\e}} + w(\e , \theta_{\e})$
and $\theta_{\e} \rightarrow \theta$ as $\e \rightarrow 0$. Assume
that $z_{\theta}= \lim_{\e} z_{\theta_{\e}} $ is nondegenerate for
the restriction of $f_{0}$ to $(T_{z_{\theta}}Z)^{\perp}$, with Morse
index equal to $m_{0}$, and that the hessian matrix $D^{2}\Gamma
(\theta )$ is positive or negative definite.
\par Then $u_{\e }$ is a nondegenerate critical point for $f_{\e }$
with Morse index equal to $m_{0}$ if $D^{2}\Gamma (\theta)$ is
positive, to $m_{0} + d$ if $D^{2}\Gamma (\theta )$ is negative. As a
consequence, the critical points of  $f_{\e }$ form a continuous
curve.
\end{Theorem}
\begin{pf}
To study the behavior of $D^{2}f_{\e}(u_{\e } )$ on $E^{+}_{\e}$ and
$E^{-}_{\e}$ we repeat the arguments of the previous theorem. As to
$E^{0}_{\e}$, we recall that the hypotheses imply $w_{\e }= O(\e^{\a
})$ (see lemma 2.2 in \cite{AB2}), so (\ref{eq:dev2}) and
(\ref{eq:dev3}) still hold and the proof goes on as in the previous
theorem
\end{pf}

\bigskip

\bigskip

\noindent We want now to apply these abstract results to our equation
(\ref{eq:gene}). In the following theorem we apply theorem
\ref{th:abstract1}. To fit hypotheses ${\bf (F_{4}) }$, ${\bf (G_{4})
}$, we have to assume $p\geq 3$. Together with the hypothesis $p <
\frac{N+2}{N-2}$, this of course implies $N \leq 3$. Notice that in 
the following theorems we will treat
curves of solutions bifurcating from $0$, or $\infty$, or bounded 
away both from $0$ and $\infty$. Recall that
we refer in our claims to the $H^{1}$-norm, and that in any case the 
$L^{\infty}$-norm is vanishing.

\bigskip

\begin{Theorem}
Let us suppose $N=1,2,3$ and $3\leq p<q< +\infty$ if $N=1,2$ while $3\leq
p<q\leq 5$ if $N=3$. Assume ${ \bf  (a_{1} ) }$,  ${ \bf  (a_{2} )
}$, ${ \bf  (b_{1} ) }$ and ${ \bf  (b_{2} ) }$.
Then we obtain a curve $(\l ,\psi_{\l } )$ of solutions of
(\ref{eq:gene}), where $\l \in (\l_{0} , 0)$, for a suitable 
$\l_{0}<0$. We have the following
behavior of $\psi_{\l}$ as $\l \rightarrow 0$:

  \begin{description}
        \item[${\bf 1.}$]   If $N=1$ and $3\leq p <5$, then $ \vert 
\vert   \psi_{\l}  \vert \vert
        \rightarrow 0$, so we
have a curve of solutions bifurcating from the origin  in $H^{1}(\re^{N})$;

      \item[${\bf 2.}$] If $N=1$ and $p=5$, or if $N=2$ and $p=3$, 
then $  \vert \vert
\psi_{\l} \vert \vert   \rightarrow c \not= 0$, so we have, in 
$H^{1}(\re^{N})$, a curve of
solutions bounded away from $0$ and $\infty$;

\item[${\bf 3.}$] If $N=1$ and $p>5$, or if $N=2$ and $p>3$, or if $N=3$ and
$3\leq p <5$, then
$  \vert \vert  \psi_{\l} \vert \vert   \rightarrow +\infty $, so we
have, in $H^{1}(\re^{N})$, a curve of solutions bifurcating from infinity.
\end{description}

\end{Theorem}
\begin{pf}
By theorem \ref{th:1} we get for equation (\ref{eq:gene1}) a family of
solutions $u_{\e} = z_{\theta_{\e}} + w( \e , \theta_{\e} )$. Using
the general devices of section 2, it is easy to see that $\theta_{\e
}$ must converge, as $\e \rightarrow 0$, to a maximum or a minimum
point of $\Gamma$. But $\Gamma$, as defined in lemma \ref{le:3}, has
the unique critical point $\theta = 0$, hence $\theta_{\e}
\rightarrow 0$ as $\e \rightarrow 0$. To apply theorem
\ref{th:abstract1} we have to show that the hypotheses ${\bf (F_{4}
)}$, ${\bf (G_{4} )}$, ${\bf (G_{5} )}$ are satisfied. The assumption
$p\geq 3$ gives that $F,G$ are $C^{4}$. As to the asymptotic
assumptions, they are easily checked with arguments similar to those
of section 3 and we leave this to the reader (notice that here $\a = 
N$). Recalling that $z_{0}$ is
a nondegenerate critical point for $f_{0}$ we apply theorem 5.1 and 
we find a curve
$(\e , u_{\e })$ of solutions of (\ref{eq:gene1}), such that
$u_{\e}(x) = z_{0}(x+\theta_{\e }) + w (\e ,\theta_{\e } ) (x)$ with 
$\theta_{\e }
\rightarrow 0$ and $ w (\e ,\theta_{\e } ) \rightarrow 0$ as $\e 
\rightarrow 0$. By the change of variables
$\l = -\e^{2}$ and $\psi_{\l} (x) =
\e^{2 /p-1}u_{\e} (\e x)$ we get a family of solutions of
(\ref{eq:gene}), which is still a curve. Noticing that $u_{\e}
\rightarrow z_{0}$ in $H^{1}$, it is easy to verify the statements on
$\lim_{\l \rightarrow 0} \vert \vert  \psi_{\l}  \vert \vert $: it is just a
computation involving a change of variables, and we leave it
to the reader.
\end{pf}

\bigskip

\noindent In the case $N=1$ it is possible to relax some hypotheses.
In the following theorem we assume $p\geq 2$ and we do not suppose
$a$ continuous and bounded. We apply theorem 5.2.

\bigskip

\begin{Theorem}\label{th:appl1}
Let us assume $N=1$ and $2\leq p<q<+\infty$. Assume $a-A \in
L^{1}(\re)$, $\int_{\re} (a(x)-A) dx \not= 0$, and suppose also ${\bf
(b_{1})}$ and ${\bf (b_{2})}$. Then we obtain a curve $(\l ,\psi_{\l
} )$ of solutions of (\ref{eq:gene}), where $\l \in (\l_{0} , 0)$, 
for a suitable $\l_{0}<0$.
We have the following behavior of $\psi_{\l}$ as $\l \rightarrow 0$:

  \begin{description}
        \item[${\bf 1.}$]  If $2\leq p <5$, then $\vert \vert 
\psi_{\l} \vert \vert \rightarrow 0$,
        so we have a
curve of solutions bifurcating from the origin in $H^{1}(\re^{N})$.

\item[${\bf 2.}$]  If $p=5$ then $  \vert \vert  \psi_{\l} \vert 
\vert   \rightarrow
c \not= 0$, so we have, in $H^{1}(\re^{N})$, a curve of solutions 
bounded away from $0$ and
$\infty$.

\item[${\bf 3.}$]    If  $p>5$ then $  \vert \vert  \psi_{\l} \vert \vert
\rightarrow +\infty $, so we have, in $H^{1}(\re^{N})$, a curve of 
solutions bifurcating
from infinity.
\end{description}

\end{Theorem}\label{th:appl2}
\begin{pf}
We want to apply theorem \ref{th:abstract2}. As $p\geq 2$, $F,G $ are
$C^{3}$. As to the asymptotic properties of $G$, in particular ${ \bf 
( G_{5}  )'}$, notice first that here
$\a =1$. We use the arguments of \cite{AB2} (in particular the proof 
of lemma 4.1, p. 1142-1143) to study
$G_{1}'$ and $G_{1}''$, while the study of $G_{2}'$ and $G_{2}''$ is 
the same as in the previous sections.
Hence we obtain a curve $(\e , u_{\e} )$
of solutions of (\ref{eq:gene1}). Arguing as in the previous theorem,
we obtain a curve $(\l ,\psi_{\l})$ of solutions of
(\ref{eq:gene}), and we get its asymptotic properties as $\l
\rightarrow 0$.
\end{pf}

\bigskip

\begin{Remark}
Theorems \ref{th:abstract1} and \ref{th:abstract2} fill a gap in the proof
of theorem 3.2 in \cite{AB2}. In that paper theorem 3.2 was used
only in theorem 1.5, to prove that a family of solutions was a curve.
Now this result is a particular case of theorem \ref{th:appl2}. A 
similar correction to theorem 3.2 of
\cite{AB2} was obtained by S. Kr\"omer in his Diplomarbeit \cite{Kr}. 
He also obtained there a
bifurcation result analogous to theorem \ref{th:1}.
\end{Remark}

\end{document}